\documentclass{article}
\usepackage[left=2.5cm,right=2.5cm,top=2.5cm,bottom=2.5cm]{geometry}
\usepackage{graphicx}
\usepackage{amsmath}
\usepackage{amssymb}
\usepackage{epstopdf}
\usepackage{hyperref}

\usepackage{color}
\usepackage{xcolor}

\newtheorem{theorem}{Theorem}
\newtheorem{lemma}[theorem]{Lemma}

\renewcommand{\d}{\mathrm{d}}

\DeclareMathOperator{\lspan}{span}
\newcommand{\R}{\mathbb{R}}

\newcommand{\tmax}{\tau_{\max}}
\newcommand{\badj}{{B^\dagger}}

\makeatother

\begin{document}
	
	\title{Nonlinear effects of instantaneous and delayed state dependence\\ in a delayed feedback loop}
	
	\author{Antony R. Humphries$^{1}$, Bernd Krauskopf$^{2}$, Stefan Ruschel$^{2}$\\ and Jan Sieber$^{3}$ \\[0.5cm] 
		$^1$ Departments of Mathematics \& Statistics, and, Physiology,\\ McGill University, Montreal, Quebec H3A 0B9, Canada\\[.2cm]
		$^2$ Department of Mathematics and Dodd-Walls Centre for Photonic and\\ Quantum Technologies 
		University of Auckland, Private Bag 92019, \\ Auckland 1142, New Zealand \\[.2cm]
		$^3$ College of Engineering, Mathematics and Physical Sciences,\\ University of Exeter, Harrison Building, Exeter EX4 4QF, United Kingdom}
	
	\date{October 2021}
	
	\maketitle


\begin{abstract}
We study a scalar, first-order delay differential equation (DDE) with instantaneous and state-dependent delayed feedback, which itself may be delayed. The state dependence introduces nonlinearity into an otherwise linear system. We investigate  the ensuing nonlinear dynamics with the case of instantaneous state dependence as our starting point. We present the bifurcation diagram in the parameter plane of the two feedback strengths showing how periodic orbits bifurcate from a curve of Hopf bifurcations and disappear along a curve where both period and amplitude grow beyond bound as the orbits become saw-tooth shaped. We then `switch on' the delay within the state-dependent feedback term, reflected by a parameter $b>0$. Our main conclusion is that the new parameter $b$ has an immediate effect: as soon as $b>0$ the bifurcation diagram for $b=0$ changes qualitatively and, specifically, the nature of the limiting saw-tooth shaped periodic orbits changes. Moreover, we show --- numerically and through center manifold analysis --- that a degeneracy at $b=1/3$ of an equilibrium with a double real eigenvalue zero leads to a further qualitative change and acts as an organizing center for the bifurcation diagram.

Our results demonstrate that state dependence in delayed feedback terms may give rise to new dynamics and, moreover, that the observed dynamics may change significantly when the state-dependent feedback depends on past states of the system. This is expected to have implications for models arising in different application contexts, such as models of human balancing and conceptual climate models of delayed action oscillator type. 
\end{abstract}

\section{Introduction}
\label{sec:intro}

A classic scenario for the potentially destabilizing effects of delays is the interplay between two feedback loops, where one is (effectively) instantaneous while the other is subject to considerable delay. Specific examples arise in numerous application areas, including balancing and control \cite{stepan}, hematopoiesis \cite{DeSouza2019}, machining \cite{Ins-Ste-00}, laser dynamics \cite{kaneshore,chapter_k} and climate modeling \cite{dijkstra,Kea-Kra-Pos-17,ks}, and they are described mathematically by delay differential equations (DDEs). We consider here the case of a scalar first-order DDE, which gives rise to the general form
\begin{equation}
\label{eq:lin}
u^{\prime}(t)=\alpha u(t)+\beta u(t-\tau)
\end{equation}
at the linear level; note that we may fix the nonzero delay at $\tau=1$ without loss of generality by choosing $\tau > 0$ as the unit of time. This linear DDE was first studied by Hayes \cite{Hayes1950} and is now used as a starting point and motivation in textbooks on DDEs such as Verduyn Lunel and Hale \cite{Hale1993} as well as Smith \cite{HalSmith2011}; see also Breda, Maset and Vermiglio \cite{Bre-Mas-Ver-14}.

Equation~\eqref{eq:lin} can be considered as a model for balancing with feedback control in the limit of low inertia (or large friction) \cite{stepan}. Here $u^{\prime}=\alpha u$ with $\alpha>0$ is the uncontrolled unstable system, where $u$ measures the deviation from the upright, and the term $\beta u(t-\tau)$ with $\beta<0$ is the stabilizing balancing force acting with a reaction delay \cite{stepan}. Similarly, conceptual models of delayed action oscillator type for the formation of global climate phenomena, such as the El Ni{\~n}o Southern Oscillation (ENSO) \cite{Sua-Sch-88,dijkstra} and the Atlantic Meridional Overturning Oscillation (AMOC) \cite{dijkstra,Rah-95} are known to feature near-instantaneous positive feedback as well as delayed negative feedback (arising from global energy transport across oceans). Such DDE models arising in applications generally also feature nonlinear terms, but their linearizations may take the form of \eqref{eq:lin}. In other words, the linear stability diagram in the $(\alpha,\beta)$-plane of system \eqref{eq:lin}, which will be discussed below in Sec.~\ref{sec:linstab}, is the starting point of the bifurcation analysis of related nonlinear DDEs.

We consider a version of \eqref{eq:lin} with state-dependent delay given by
\begin{equation}
\label{eq:sdDDE_a_eta}
u^{\prime}(t)=\alpha u(t)+\beta u(t-a-\eta u(t-b)) \mbox{.}
\end{equation}
Here the delay $\tau = a+\eta u(t-b)$ depends linearly on the state $u$ with strength $\eta$. We are interested in system~\eqref{eq:sdDDE_a_eta} with a nonzero delay $a > 0$ in the presence of state dependence, that is, $\eta > 0$. Rescaling $u$ and time $t$ allows us to set $a=\eta=1$ without loss of generality, yielding the state-dependent DDE
\begin{equation}
\label{eq:sdDDE}
u^{\prime}(t)=\alpha u(t)+\beta u(t-1-u(t-b))
\end{equation}
as our central object of study, where $\alpha$, $\beta$ and $b$ are the bifurcation parameters, and we restrict ourselves to the physically relevant parameter range $\alpha+\beta<0.$

Of particular interest from a purely mathematical as well as from an application point of view is the fact that the state dependence of the delay introduces nonlinearity.
This means that the system~\eqref{eq:sdDDE} may feature dynamics beyond its non-state-dependent cousin, Eq.~\eqref{eq:lin}. It has been shown \cite{Calleja2017,HumphriesJDDE16,Hum-Dem-Mag-Uph-12,JMPRNP94} that state dependence alone is capable of generating a wealth of dynamical phenomena, including resonant and multi-frequency behavior, in a scalar DDE with two or more state-dependent feedback terms. The special case of system~\eqref{eq:sdDDE} with $b=0$ was introduced by Mallet-Paret and Nussbaum. In \cite{JMPRN11} they investigated the existence and saw-tooth shaped form of the slowly oscillating periodic solutions of a singularly perturbed version of \eqref{eq:sdDDE_a_eta} with $\alpha\ll0$, $\beta\ll0$ and $b=0$, and used the same equation as an illustrative example of more general problems in \cite{JMPRNI,JMPRNII,JMPRNIII} (see also \cite{Kozyreff2014,Hum-Mag-2020}). Moreover, Magpantay and Humphries \cite{Mag-Hum-20} considered the system~\eqref{eq:sdDDE} with $b=0$ and provided sufficient conditions for the existence of periodic solutions, thus, identifying regions in the $(\alpha,\beta)$-plane where they may be found.

The novel aspect of the system~\eqref{eq:sdDDE} is the parameter $b$, which is included so that the delay may depend not only on the current state $u(t)$ but instead on the state $u(t-b)$ some given time $b\geq 0$ ago. This is motivated by the observation in \cite{Keane2019} that a delayed state dependent feedback term arises in the modeling of ENSO when one takes into account the vertical heat transport in the ocean relating the thermocline depth to the sea-surface temperature. 
Apart from \cite{Golubenets21,Keane2019}, there has been very little consideration to date of delay differential equation models with state-in-the-past dependent delays. In particular, the dynamics and bifurcations of the specific system~\eqref{eq:sdDDE} have not previously been studied for $b>0$. We mention here that Kennedy \cite{kennedy2019poincare} showed for a more general class of delayed state-dependent but negative feedback that the eventual dynamics are essentially two-dimensional. This result applies to Eq.~\eqref{eq:sdDDE} for the case $\alpha<0$ and $\beta<0$ of negative feedback; thus extending to $b>0$ earlier results by Krisztin and Arino \cite{Krisztin2001} on instantaneous state dependence.

Our goal here is to present the bifurcation analysis of system~\eqref{eq:sdDDE}, which we achieve by determining its bifurcation diagram in the $(\alpha,\beta)$-plane for different values of the delay $b$ of the state-dependent feedback. To this end, we combine analytical techniques and state-of-the-art numerical continuation for DDEs with state-dependence, which have both become available quite recently.

The classic theory of DDEs with constant delays can be found in textbooks, such as \cite{Diekmann1995,Hale1993,HalSmith2011,stepan}.  The general framework for treating DDEs as dynamical systems is to consider the current \emph{history segment} $u_t:[-\tmax,0]\to\R$, defined by $u_t(s)=u(t+s)$ as the state of the dynamical system with an extension rule; for system~\eqref{eq:sdDDE} it takes the form
\begin{align}\label{eq:gendde}
  u'(t)=f(u_t), \quad \mbox{ with } \quad f(\phi)=\alpha \phi(0)+\beta \phi(-1-\phi(-b))\mbox{,}
\end{align}
where $f:C([-\tmax,0],\R)\to\R $ and $C^0 = C([-\tmax,0],\R)$ denotes the space of continuous functions. 

For functionals $f$ that are differentiable at least $k$ times on $C^0$ the DDE $u'(t)=f(u_t)$  generates a $k$-times differentiable (continuous-time) dynamical system or flow given by the solution map $\Phi^t: u_0\mapsto u_t$. However, this does not apply to DDEs with state dependency. The functional $f$ in \eqref{eq:gendde} is not differentiable on $C^0$ due to the presence of the composition term; in fact, $f$ is not even locally Lipschitz continuous. As a result, even if we assume that $u$ stays bounded away from $-1$ and $+\infty$ such that the delay $1+u(t-b)$ of system \eqref{eq:sdDDE} stays in a finite interval $[\tau_{\min},\tmax]$ with $0<\tau_{\min}<\tmax<\infty$, the flow does not possess the regularity properties one would expect given that all terms of the right-hand side of \eqref{eq:gendde} are nested linear functions. Rather, functionals $f$ of state-dependent DDEs satisfy a weaker condition \cite{Hartung2006}, sometimes called \emph{mild differentiability} \cite{Hartung2006,Sieber12,Sieber2017}. With the help of this concept, DDEs of type \eqref{eq:gendde} were proven in \cite{Walther03} to be $C^1$-regular (that is, once differentiable) dynamical systems on the phase space (manifold) of compatible initial conditions $C^1_\mathrm{comp}=\{u\in C^1:u'(0)=f(u)\}$. This implies basic results, such as the principle of linearization along trajectories and the existence of $C^1$ local center-unstable manifolds \cite{Qesmi2009,Stu11}. In fact, in \cite{Sieber12} it is shown that periodic orbits and their bifurcations are given by roots of finite-dimensional algebraic systems of equations that have $C^k$ smooth coefficients. Furthermore, center manifold expansions and normal forms at equilibria can still be computed, and their predictions concerning families of periodic orbits branching off from the equilibrium still hold \cite{Calleja2017,Sieber2017}.

Similarly, continuation tools for DDEs with constant delays have been available for some time, as implemented in the program {\tt knut} \cite{Szalai}
and the MATLAB package DDE-BIFTOOL \cite{NewDDEBiftool}, the latter of which we use in this work. Their capabilities include the continuation of steady states, periodic orbits and their codimension-one bifurcations; see \cite{sr,CISMchapter} for background information. The present version of DDE-BIFTOOL \cite{NewDDEBiftool} is able to compute center manifold expansions of bifurcations of steady states of codimension up to two \cite{BJK20}. Importantly for this work, the formulation of the DDE in the package DDE-BIFTOOL allows for state-dependent delays, so that the same suite of capabilities is available for state-dependent DDEs \cite{NewDDEBiftool,CISMchapter}.

With these analytical and numerical tools to hand, the starting point for our study of system~\eqref{eq:sdDDE} is the case $b=0$ of instantaneous state dependence. We extend this earlier work by presenting the bounding curves of the regions of existence of periodic orbits that bifurcate from a Hopf bifurcation and then quickly become saw-tooth shaped. The existence of this type of periodic orbits indicates a definite time-scale separation along the orbit despite the absence of small or large parameters. This observation has been a motivation behind the interest in state-dependent delays from the mathematical perspective, and results on existence, regularity and monotonicity properties of such saw-tooth shaped periodic orbits can be found in \cite{JMPRNI,JMPRNII,JMPRNIII,JMPRN11}. We show how the periodic orbits of system~\eqref{eq:sdDDE} reach different limits, associated with different slopes of their saw-tooth limiting shape, as their period and amplitude go to infinity. These limits correspond to singularities when the state-dependent delay becomes advanced, or equivalently, when the periodic orbit attains its minimum value at a certain threshold value. The respective loci of these singularities are boundary curves of the region of existence of periodic orbits in the $(\alpha,\beta)$-plane for $b=0$.

We next show how the bifurcation diagram in the $(\alpha,\beta)$-plane changes when $b$ is increased from $0$. Our main conclusion is that $b>0$ immediately influences the large-period periodic behavior of the saw-tooth orbits for $b=0$ and, hence, the entire bifurcation diagram. We identify as an organizing center for this type of dynamics a point we refer to as \textbf{DZ} that corresponds to an infinitely degenerate case of an equilibrium with a double real eigenvalue zero. We determine a finite-order expansion of the two-dimensional local center manifold to analyze the dynamics near the special point \textbf{DZ}. Symbolic algebra drivers for the expansion of center manifolds in DDEs with state-dependent discrete delays and the resulting expressions up to order $5$ near \textbf{DZ} are provided as supplementary material. Even though the smoothness of local center manifolds in DDEs with state-dependent delays is not completely settled \cite{Kri-06}, on the two-dimensional center manifold the Poincar{\'e}-Bendixson theorem ensures that knowledge about periodic orbits is sufficient to capture all of the dynamics locally near \textbf{DZ} as a function of all three parameters $\alpha$, $\beta$ and $b$ of system~\eqref{eq:sdDDE}.

This paper is organized as follows. Section~\ref{sec:bif-analysis-0} is dedicated to the bifurcation analysis of the $(\alpha,\beta)$-plane for $b=0$. Here, our starting point in Sec.~\ref{sec:linstab} is the linear stability analysis of Eq.~\eqref{eq:lin}, which is followed in Sec.~\ref{sec:nonlinstab} by the presentation of the bifurcation diagram of the nonlinear system~\eqref{eq:sdDDE}; the latter features the emergence of saw-tooth periodic orbits at the Hopf-bifurcation curve, whose limiting behavior is discussed in Sec.~\ref{sec:saw-tooth-limits}. How the bifurcation diagram in the $(\alpha,\beta)$-plane changes for $b > 0$ is the subject of Sec.~\ref{sec:bif-analysis}; it involves the continuation of periodic orbits up to large periods to determine how they bifurcate. Section~\ref{sec:degDZ} presents the center manifold analysis of the degenerate equilibrium. In Section~\ref{sec:conclusions}, we conclude and point out some directions for future research. Appendix~\ref{app:ex} provides technical details regarding the computation of the  expansions of the ODE on the center manifold near the point \textbf{DZ}, and Appendix~\ref{app:supp} points out that a general implementation and the specific symbolic output is available as supplementary material.

\section{Bifurcation analysis of instantaneous state-dependency, $\boldsymbol{b=0}$}
\label{sec:bif-analysis-0}
We first consider the classical situation that the state-dependent term is not delayed itself, that is, $b=0$ in Eq.~\eqref{eq:sdDDE}. The starting point is the linear stability analysis of the equilibrium $u\equiv0$ of Eq.~\eqref{eq:lin}, which is indeed also valid for $b>0$. We then consider nonlinear stability, that is, the bifurcation diagram of system~\eqref{eq:sdDDE} for $b=0$ in the $(\alpha,\beta)$-plane. This involves the study of periodic orbits bifurcating from a Hopf bifurcation and their disappearance along curves where their period goes to infinity.

\subsection{Linear stability analysis}
\label{sec:linstab}

The linearization of Eq.~\eqref{eq:sdDDE} about the equilibrium $u\equiv0$ was first studied in \cite{Hayes1950} and is now a standard example in DDE textbooks \cite{Bre-Mas-Ver-14,Hale1993,HalSmith2011,stepan}. It is given by 
\begin{equation}
\label{eq:lin_tau1}
x^{\prime}(t)=\alpha x(t)+\beta x(t-1),
\end{equation}
does not depend on the delay $b$ of the state-dependent term and is indeed Eq.~\eqref{eq:lin} with $\tau=1$. The characteristic equation of Eq.~\eqref{eq:lin_tau1} for the eigenvalues of $u\equiv0$ is
\begin{equation}
\label{eq:char}
\lambda-\alpha-\beta e^{-\lambda}=0.
\end{equation}
It follows immediately that there is a real eigenvalue $\lambda=0$ along the line $\beta = - \alpha$ in the $(\alpha,\beta)$-plane, which we denote \textbf{Z}. Moreover, there may be complex conjugate pairs $\lambda$, $\overline{\lambda}$ of eigenvalues with zero real part, and this happens along infinitely many curves in parameter space, which we refer to as $H_n$ where $n = 0, 1, 2, \ldots$; see, for example, \cite{HalSmith2011,Bre-Mas-Ver-14} where more details can be found. Along the curve \textbf{H} $= H_0$ given by
\begin{equation}
\label{eq:C0}
\mbox{\textbf{H}}=\Bigl\{(\alpha,\beta) =(\theta\cot\theta, -\theta\csc\theta), \mbox{ where } \theta\in[0,\pi)\Bigr\}
\end{equation}
the equilibrium $u\equiv0$ loses its stability. The curve \textbf{H} starts on the curve \textbf{Z} at the point $(\alpha,\beta)=(1,-1)$ for $\theta=0$, crosses the $\beta$-axis of the $(\alpha,\beta)$-plane at $\beta=-\pi/2$ and approaches the diagonal as both $\alpha$ and $\beta$ go to $-\infty$ as $\theta \to \pi$. At the point $(\alpha,\beta)=(1,-1)$ the equilibrium $u\equiv0$ has a double real zero eigenvalue ($\lambda=0$ is a double root), and hence we refer to this point as \textbf{DZ}.

\begin{figure}[t!]
	\centering
	\includegraphics[width=0.76\textwidth]{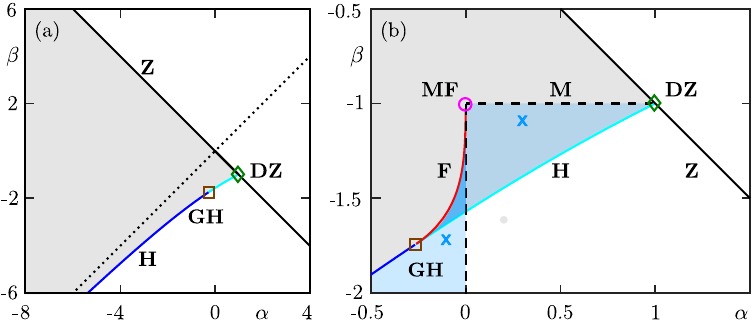}
	\caption{Bifurcation diagram of system~\eqref{eq:sdDDE} with $b=0$ in the $(\alpha,\beta)$-plane. Panel~(a) shows the stability diagram of the equilibrium $u\equiv0$ with the curves \textbf{Z} (black) of a simple zero eigenvalue and \textbf{H} of a pair of complex conjugate eigenvalues with zero real part, which meet at the \textbf{DZ} (green diamond) of double zero eigenvalues and bound the region where \eqref{eq:sdDDE} is linearly stable (shaded); also shown for reference is the line $\alpha=\beta$ (black dotted). The associated Hopf bifurcation of \eqref{eq:sdDDE} along the curve \textbf{H} has a generalized Hopf bifurcation point \textbf{GH} (brown square), where its criticality changes from supercritical (blue) to subcritical (cyan). Panel~(b) is an enlargement that also shows a curve of folds of periodic orbits \textbf{F} (red), a curve \textbf{M} (black dashed) of periodic orbits with a point that attains the minimum value $-1$, which has a corner at \textbf{MF} (magenta circle) where there is a fold periodic orbit with minimum value $-1$. The blue shaded regions indicate the (co)existence of periodic orbits, and the blue crosses mark the locations of the periodic orbits shown in Fig.~\ref{fig:periodic}.}
\label{fig:stabdiag_b0}
\end{figure}

These curves and the stability region are shown in Fig.~\ref{fig:stabdiag_b0}(a).
The equilibrium $u\equiv0$ of the linear system~\eqref{eq:lin} is exponentially stable in the shaded stability region in the $(\alpha,\beta)$-plane, which is bounded by the curves \textbf{Z} and \textbf{H}; this region is open to the left and its tip is the point \textbf{DZ}.
This type of stability region is common in dynamical systems with delay. Note that the lower boundary curve \textbf{H} of the stability region of the equilibrium is due to the nonzero delay, and it is of particular interest when studying delay-induced oscillatory instability in engineering \cite{stepan} and science \cite{Sua-Sch-88}.

The additional curves $H_n$ along which one finds further pairs of complex conjugate eigenvalues with zero real parts have the same form as \textbf{H} in \eqref{eq:C0}, except that $\theta\in(n\pi,(n+1)\pi)$ for $n>0$. We do not consider them here because none of the curves $H_n$ form part of the stability boundary. They are interleaved in the $(\alpha,\beta)$-plane, with no intersections with one another, and all lie outside the part of the $(\alpha,\beta)$-plane shown in Fig.~\ref{fig:stabdiag_b0}(a); see, for example, \cite{Bre-Mas-Ver-14,HalSmith2011} for illustrations of the curves $H_n$.

\subsection{Nonlinear stability and bifurcating periodic orbits}
\label{sec:nonlinstab}

In the presence of nonlinearity, global stability is no longer guaranteed and the linearization \eqref{eq:lin} only provides the information that the equilibrium is locally stable in its stability region. For system~\eqref{eq:sdDDE}, where the nonlinearity is generated by the state dependence, it has been shown in \cite{Mag-Hum-20} that the equilibrium $u\equiv0$ is globally asymptotically stable in the cone $|\beta|+\alpha<0$ bounded to the right by the lines \textbf{Z} and $\alpha=\beta$ --- but not in the whole linear stability region due to a co-existence of the stable equilibrium with a stable periodic orbit.

We now review these results and then present the bifurcation diagram in the $(\alpha,\beta)$-plane of the nonlinear system~\eqref{eq:sdDDE}, including regions of (co)existence of stable and unstable periodic solutions. First of all, the stability boundaries \textbf{Z} and \textbf{H} become bifurcation curves. Namely,  along the curve \textbf{Z} there is an infinitely degenerate transcritical bifurcation. This bifurcation is degenerate because for $\alpha=-\beta$ system~\eqref{eq:sdDDE} has a one-parameter family of equilibria, $x\equiv x_{\mathrm{eq}}$ with arbitrary $x_{\mathrm{eq}}$. Below the line $\alpha=\beta$ one finds additional interesting dynamics. As was already mentioned, along the curve \textbf{H} one encounters a Hopf bifurcation. DDE-BIFTOOL \cite{NewDDEBiftool} computes the coefficient of the third-order expansion for Hopf bifurcations, which allows one to determine the first Lyapunov coefficient and, therefore, the criticality of the Hopf bifurcation along \textbf{H}. As shown in Fig.~\ref{fig:stabdiag_b0}(a) and in the enlargement in Fig.~\ref{fig:stabdiag_b0}(b), the criticality changes at a generalized Hopf bifurcation point \textbf{GH} located at $(\alpha_{GH},\beta_{GH})\approx(-0.63,-1.99)$. In between the points \textbf{GH} and \textbf{DZ} the Hopf bifurcation along \textbf{H} is subcritical and the unstable periodic orbit that bifurcates coexists with the stable equilibrium $u\equiv0$ in the region above the curve \textbf{H}; in this case the equilibrium cannot be globally asymptotically stable. Otherwise, for $\alpha<\alpha_{GH}$, the Hopf bifurcation along \textbf{H} is supercritical and a stable periodic orbit bifurcates, which coexists with the unstable equilibrium $u\equiv0$ in the region below \textbf{H}.

\begin{figure}[t!]
	\centering	\includegraphics[width=0.76\textwidth]{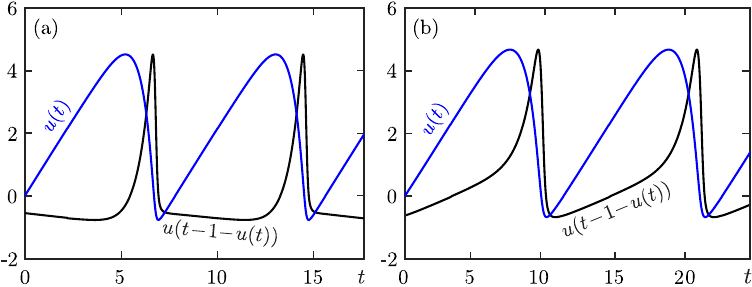}
	\caption{Periodic solutions of system~\eqref{eq:sdDDE} for $b=0$ shown as time series of $u(t)$ and of $u(t-1-u(t))$; panel~(a) shows the stable periodic orbit for $\alpha=-0.1$ and $\beta=-1.7$, and panel~(b) the unstable orbit for $\alpha=0.3$ and $\beta=-1.1$.
These parameter values are denoted by blue crosses in Fig.~\ref{fig:stabdiag_b0}.}
\label{fig:periodic}
\end{figure}

Figure~\ref{fig:periodic} shows time series of a stable and of an unstable periodic orbit, namely those that exist at the two crosses in the $(\alpha,\beta)$-plane of Fig.~\ref{fig:stabdiag_b0}(b), where the existence of periodic orbits is indicated by additional shading. Figure~\ref{fig:periodic} shows not only the time series of $u(t)$ but also that of the delayed variable $u(t-1 - u(t))$. Observe that this type of oscillation follows closely the straight line $u(t)\sim k_1t$ with a positive slope $k_1$ up to a large value of $u$, before a quite sharp transition back to $u\approx -1$. The latter is accompanied by a spike and then reset of $u(t-1 - u(t))$, where the delay $1+u(t)$ almost vanishes when $u\approx -1$. The difference between the stable periodic orbit in panel~(a) and the unstable periodic orbit in panel~(b) of Fig.~\ref{fig:periodic} lies in the subsequent slope of $u(t-1 - u(t))$ after the spike.

Going beyond the results in~\cite{Mag-Hum-20}, we now perform a continuation of the stable periodic orbit bifurcating from \textbf{H}. It shows that, as $\alpha$ is increased for fixed $\beta$, the unstable periodic orbit disappears along the half-line given by $\alpha = 0$ with $\beta < -1$. Similarly, as $\beta$ is increased for fixed $\alpha$, the unstable periodic orbit bifurcating from \textbf{H} disappears along the line segment given by $\beta = -1$ with $0 < \alpha < 1$. Together these two curves form the locus \textbf{M} along which the respective periodic orbit has a point that attains the minimum value $-1$. As is discussed in the next section, the stable and unstable periodic orbits, which already appear to be quite close to a saw-tooth shape in Fig.~\ref{fig:periodic}, reach different limiting shapes as these two branches of \textbf{M} are approached, which implies that the locus \textbf{M} is a boundary for the existence of the respective periodic orbit.

We remark that the limiting periodic orbits close to \textbf{M} are reminiscent of the saw-tooth periodic orbits studied by Mallet-Paret and Nussbaum \cite{JMPRN11}. To be more specific, those authors proved rigorously the existence of slowly oscillating periodic orbits with saw-tooth shape of the state-dependent delay differential equation $\varepsilon u^{\prime}(t)= -u(t)-ku(t-1-u(t))$ when $\varepsilon$ is sufficiently small. We note that Eq.~\eqref{eq:sdDDE} for $b=0$ can be brought into this form by setting $\alpha=-1/\epsilon$ and $\beta=-k/\varepsilon$, subject to the restriction $\alpha<0$ and $\beta<0$.
This apparent connection motivates us to study in Sec.~\ref{sec:saw-tooth-limits} the scaling behavior of periodic orbits near \textbf{M} in more detail, that is for both $\alpha=0$ and $\beta\leq-1$ as well as $\beta=-1$ and $0\leq\alpha<1$. A rigorous analysis in the spirit of Mallet-Paret and Nussbaum \cite{JMPRN11} is beyond the scope of our work here, but could present an interesting and challenging topic for future research.

The boundary of the region with periodic orbits also contains the curve \textbf{F} of fold periodic orbits, which emerges from the generalized Hopf bifurcation point \textbf{GH} (as predicted by its normal form \cite{Kuz-04-Book}). Continuation of \textbf{F} shows that this curve ends at the point labeled \textbf{MF} at $(\alpha_{MF},\beta_{MF})\approx(0,-1)$, which is the corner point of the locus \textbf{M}. In the triangular region bounded by the curves \textbf{H}, \textbf{F} and \textbf{M} in Fig.~\ref{fig:stabdiag_b0}(b) the stable periodic orbit coexists with the stable equilibrium $u\equiv0$, as well as the unstable periodic orbit; the two periodic orbits bifurcate and disappear as this region is exited through the fold curve \textbf{F}. We remark that the curve \textbf{F}, but not the curve \textbf{M}, was found previously in \cite{Mag-Hum-20}.

\subsection{Two different saw-tooth limits and the locus \textbf{M}}
\label{sec:saw-tooth-limits}

We now consider the limits of the stable and unstable periodic orbits in Fig.~\ref{fig:periodic} as either $\alpha$ increases to $0$ for fixed $\beta < -1$, or  $\beta $ increases to $-1$ for fixed $0 \leq \alpha < 1$, respectively. To this end, we perform one-parameter continuation of these periodic orbits in the respective parameter, starting at the Hopf bifurcation curve \textbf{H}, until either a large period $T$ is reached or the minimum of periodic orbit comes very close to the value $-1$ (so that the state-dependent delay becomes zero). These and related quantities are monitored, which allows us to characterize the nature of the respective limiting process. In particular, we confirm in this way the nature of the locus \textbf{M} in Fig.~\ref{fig:stabdiag_b0}(b).

\subsubsection{Limit of the stable periodic orbit as $\alpha$ increases to $0$}
\label{sec:limit_alpha0}

\begin{figure}[t!]
	\centering
	\includegraphics[width=0.76\textwidth]{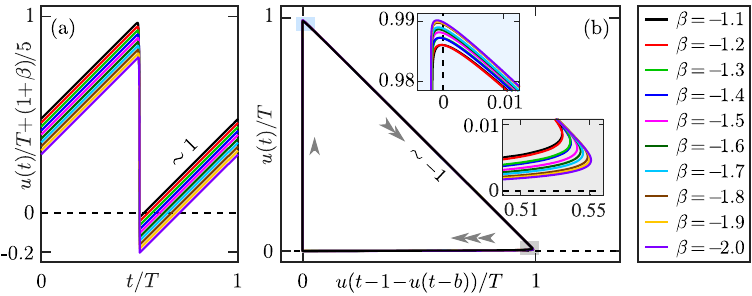}
	\caption{Periodic orbits of system~\eqref{eq:sdDDE} for $b=0$ as $\alpha$ increases to $0$ for different fixed values of $\beta<-1$ as shown. Panel~(a) shows their asymptotic profiles in a waterfall plot synchronized to their fast segments, and panel~(b) shows them in the $(u(t),u(t-1-u(t)))$-plane with enlargements near the highlighted two corners; all panels are scaled by the period $T$, the periodic orbits have been continued in to $T=5\cdot10^{2}$, and the arrows in panel~(b) indicate the time-scale separation along the orbit.}
\label{fig:limit_alpha0}
\end{figure}

Figure~\ref{fig:limit_alpha0} shows the ten periodic orbits of system~\eqref{eq:sdDDE} for $b=0$ and $\beta = -1.1$ down to $\beta = -2.0$ (see the legend) for which $T=5\cdot10^{2}$ has been reached as $\alpha$ approaches $0$ from below. Panel~(a) shows the respective time profiles with a small offset in a waterfall plot; here both $t$ and $u$ have been scaled by $T$ and the profiles have been synchronized so that their fast segments are at the (rescaled) time instance $t/T = 0.5$. This illustrates that for all values of $\beta$ asymptotically, as $\alpha$ increases towards $0$, the same limiting saw-tooth shape is reached, consisting of linear solution segments with slope $1$ that are periodically reset to negative values. Note that the time interval during which the switch occurs is small. This limiting behavior is illustrated further in Fig.~\ref{fig:limit_alpha0}(b), where these periodic orbits are shown in projection onto the $(u(t-1-u(t))/T,u(t)/T)$-plane. All ten periodic orbits practically coincide in this representation and effectively follow one and the same (rescaled) isosceles right triangle, where gray arrows indicate the time-scale separation during the different stages of the saw tooth orbit. The orbit mainly consists of a linearly increasing segment with slope $1$ of $u(t)$ from $0$ to $1$, while $u(t-1-u(t))$ is practically $-1$ (close to zero when rescaled by the period). This is followed by the fast reset to negative values of $u(t)$ when $u(t)$ is close to its maximum value, which is approximately given by $(1+\beta)/\alpha$ and accompanied by $u(t-1-u(t))$ changing its sign; see the upper inset of Fig.~\ref{fig:limit_alpha0}(b). During the reset there is a fast segment where $u(t)$ first decreases linearly with respect to $u(t-1-u(t))$ all the way to $0$, while $u(t-1-u(t))$ increases to its (rescaled) maximum $1$. In particular, $u(t-1-u(t))$ catches up with $u(t)$ and its dynamics is the fastest when it is close to its maximum value. During this very fast stage, $u(t-1-u(t))$ catches up with $u(t)$ for the second time, while $u(t)$ is still practically $0$ and the cycle repeats.

\begin{figure}[t!]
	\centering	\includegraphics[width=0.76\textwidth]{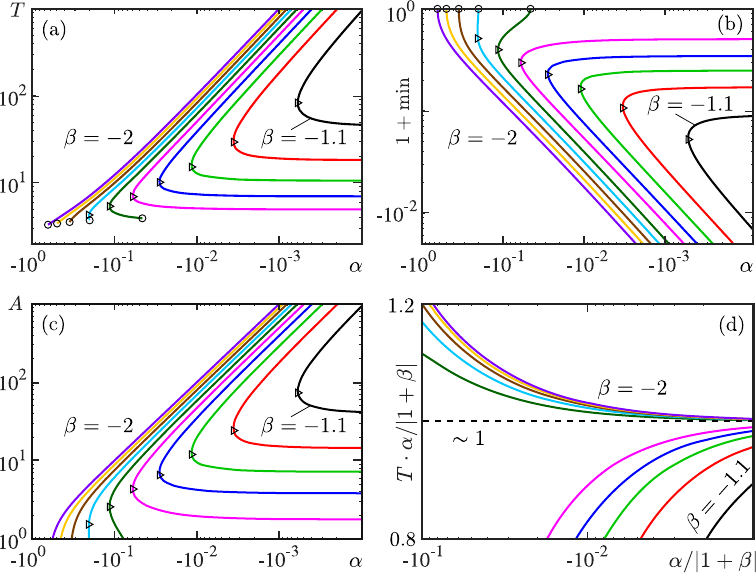}	
	\caption{Continuation of eventually stable periodic orbits of system~\eqref{eq:sdDDE} for $b=0$ in $\alpha$ for different fixed values of $\beta<-1$; compare with Fig.~\ref{fig:limit_alpha0}. Illustrated on a logarithmic scale for $\alpha$ are the asymptotic behavior of the period in panel~(a), the minimum in panel~ (b), the amplitude in panel~ (c), and the rescaled period $T \cdot \alpha / (1+\beta)$ in panel~ (d); circles indicate Hopf bifurcations and triangles fold points.}
	\label{fig:continue_alpha0}
\end{figure}

While all ten periodic orbits in Fig.~\ref{fig:limit_alpha0}(b) are extremely close to each other, small differences between the periodic orbits can be distinguished when one enlarges the turns near the corners. The two inset panels show that these turns are sharper the smaller $\beta$ is, which is because we show periodic orbits with a constant period. Namely, $T\sim (1+\beta)/\alpha$ meaning that $1+\beta$ is inversely proportional to the distance of $\alpha$ to zero (that is, to the locus \textbf{M}), which determines the interval of reset and, hence, the sharpness of the turns.

Figure~\ref{fig:continue_alpha0} provides additional information about the periodic orbits for fixed $\beta$, which are born at the Hopf bifurcation curve \textbf{H} and eventually approach the saw-tooth limit along the vertical part of the locus \textbf{M} as $\alpha$ is varied. More specifically, shown are the continuations of the ten periodic orbits for the fixed values of $\beta<-1$ from Fig.~\ref{fig:limit_alpha0}. They are represented in Fig.~\ref{fig:continue_alpha0} on a doubly-logarithmic scale by the period in panel~(a), by the $u$-minimum (shifted up by $1$) along the profile in panel~(b), and by the amplitude in panel~(c), shown as functions of $\alpha$. For $\beta=-1.1$ down to $\beta=-1.7$, the periodic orbit is born along the subcritical part of the curve \textbf{H}; hence, $\alpha$ initially decreases until the fold point on the curve \textbf{F} is reached and $\alpha$, which is now negative, increases as the curve \textbf{M} is approached; compare with Fig.~\ref{fig:stabdiag_b0}. Note that the periodic orbit is stable past the fold, as are the periodic orbits for $\beta=-1.8$ down to $\beta=-2$ that bifurcate from the supercritical part of the Hopf bifurcation curve \textbf{H} in the direction of increasing (and also negative) $\alpha$. As $\alpha$ increases to $0$ and the curve \textbf{M} is approached, period and amplitude of any of these stable periodic orbits grow with decreasing distance of $\alpha$ from $0$; at the same time the minimum along the periodic orbit approaches $-1$. Notice in Fig.~\ref{fig:continue_alpha0}(a)--(c) that all these curves have the same slope when  $\alpha$ is negative and sufficiently close to $0$. Panel~(d) illustrates that $T$ scales as $|1+\beta|/\alpha$ as $\alpha$ approaches $0$, shown as the rescaled period $T \cdot \alpha/|1+\beta|$ versus $\alpha/|1+\beta|$ on a linear scale. This scaling also holds for the amplitude, due to the linear relationship between period and amplitude shown in panels~(a) and~(c) as $\alpha$ increases to $0$.

\subsubsection{Limit of the unstable periodic orbit as $\beta$ increases to $-1$}
\label{sec:limit_betamin1}

\begin{figure}[t]
	\centering
	\includegraphics[width=0.76\textwidth]{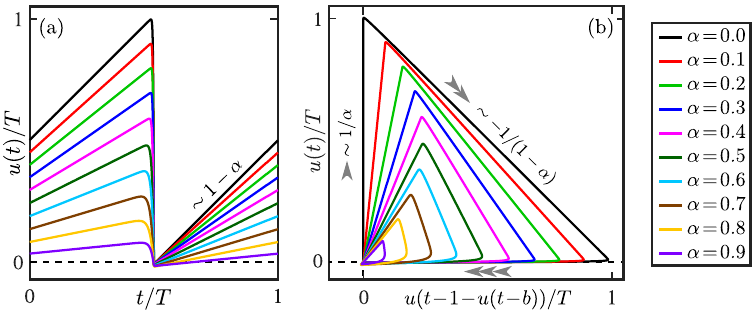}
	\caption{Periodic orbits of system~\eqref{eq:sdDDE} for $b=0$ as $\beta$ increases to $-1$ for different fixed values of $0\leq\alpha<1$, shown as in Fig.~\ref{fig:limit_alpha0}; here the periodic orbits for $0\leq\alpha\leq0.2$ have been continued to $T=5\cdot10^{2}$ and those for $0.3\leq\alpha\leq0.9$ to $\beta=-1-10^{-5}$.}
\label{fig:limit_betamin1}
\end{figure}

We now consider the behavior of the unstable periodic orbit of system~\eqref{eq:sdDDE} for $b=0$ as the horizontal segment of the locus \textbf{M} is approached as $\beta$ increases to $-1$ for different fixed values of $\alpha$ with $0 < \alpha < 1$. Figure~\ref{fig:limit_betamin1} illustrates ten computed periodic orbits near \textbf{M} (as determined by either a high period or a small value of $\beta +1$) for $\alpha = 0$ up to $\alpha = -0.9$ (see the legend). As in the previous section, panel~(a) shows the respective time profiles, where both time and $u(t)$ have been scaled by $T$ and the fast segments are at the (rescaled) time instance $t/T = 0.5$. When $\alpha=0$, which corresponds to the corner point $(0,-1)$ of the locus \textbf{M} in the $(\alpha,\beta)$-plane, the slope is $1$ as was the case in the previous section. However,  in contrast to the limit we considered in Sec.~\ref{sec:limit_alpha0}, the slope of the limiting saw-tooth shape along the horizontal part of \textbf{M} is now $1-\alpha$ and so depends on $\alpha$, which is why no offset is required to distinguish these ten orbits in Fig.~\ref{fig:limit_betamin1}(a). In particular, the limiting slope as $\beta \to -1$ approaches zero as also $\alpha\to1$, that is the end point \textbf{DZ} of the locus \textbf{M}.

The projection onto the $(u(t-1-u(t))/T,u(t)/T)$-plane in Fig.~\ref{fig:limit_betamin1}(b) shows that, owing to the different slopes involved, the ten unstable periodic orbits do not coincide. We observe that, as $\alpha$ increases, the maximum value along the respective periodic orbit decreases and scales as $A\sim(1-\alpha)T$ with period; as a result, the curve of maxima in this projection is approximately given by $(\alpha(1-\alpha),1-\alpha)$. Again, gray arrows indicate the time-scale separation during the different stages of the saw tooth orbit. During the slow part, $u(t-1-u(t))$ increases linearly with slope $\alpha(1-\alpha)$ and we find the relationship $(1+u(t))/(1-u(t-1-u(t)))\approx\alpha$ between the instantaneous and the delayed state. The fast reset to negative values of $u(t)$ occurs during a time interval of order $1$ compared to the large period when $u(t)$ is close to its maximum value; now $u(t)$ decreases linearly while $u(t-1-u(t))$ increases fast, in such a way that it gives rise to the slope $-1/(1-\alpha)$ in Fig.~\ref{fig:limit_betamin1}(b). Again, it is during the switch that $u(t-1-u(t))$ catches up with $u(t)$. The dynamics of $u(t-1-u(t))$ is the fastest when $u(t-1-u(t))$ reaches its maximum value. While $u(t)$ is still small $u(t-1-u(t))$ catches up with $u(t)$ for the second time and the cycle repeats.

Figure~\ref{fig:continue_betamin1} shows the scaling behavior of period, minimum and amplitude of the ten unstable periodic orbits as a function of $\beta$ (offset by $1$ again to enable a doubly logarithmic scale), when they are continued from the Hopf bifurcation curve \textbf{H} as $\beta$ increases toward $-1$. Note that all of these periodic orbits are unstable, since they bifurcate from the unstable part of \textbf{H} in the direction of increasing $\beta$; compare with Fig.~\ref{fig:stabdiag_b0}. As panels~(a)--(c) show, in the process, period and amplitude again grow beyond bound. However, this growth depends on the value of $\alpha$: the closer $\alpha$ is to $1$, the slower period and amplitude increase as a function of the distance of $\beta$ from $-1$; similarly, the decrease of the minimum along the periodic orbit approaches $-1$ slower the larger $\alpha$ is. This slowing growth with $\alpha$ is the reason why the periodic orbits shown in Fig.~\ref{fig:limit_betamin1} could be continued up to fixed $T=5\cdot10^{2}$ only for $0\leq\alpha\leq0.2$; those for $0.3\leq\alpha\leq0.9$ have been continued to fixed $\beta=-1-10^{-5}$ instead. Fig.~\ref{fig:continue_alpha0}(d) shows that period and amplitude indeed scale as $A\sim(1-\alpha)T$ in the limit $\beta\to-1$.

\begin{figure}[t!]
	\centering	\includegraphics[width=0.76\textwidth]{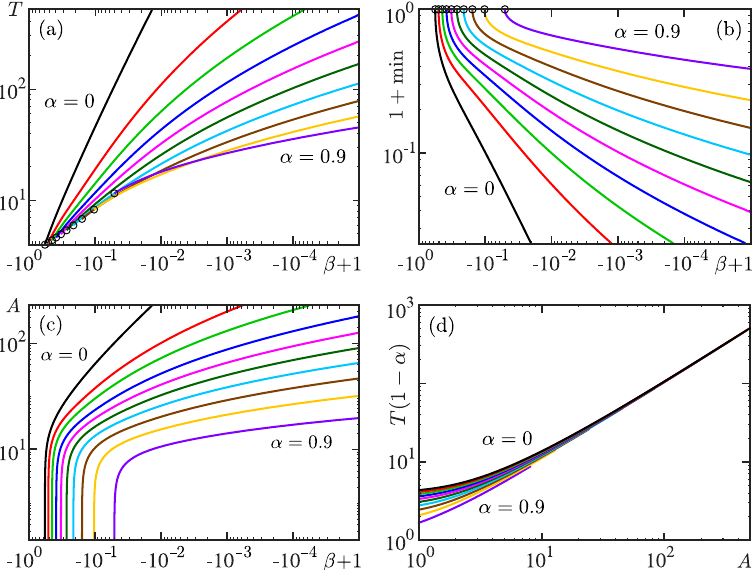}
	\caption{Continuation of unstable periodic orbits of system~\eqref{eq:sdDDE} for $b=0$ in $\beta$ for different fixed values of $\alpha$; compare with Fig.~\ref{fig:limit_betamin1}. Illustrated on a logarithmic scale for $\alpha$ are the asymptotic behavior of the period in panel~(a), the minimum in panel~(b), the amplitude in panel~(c), and the rescaled period $T(1-\alpha)$ in panel~(d); circles indicate Hopf bifurcations.}
	\label{fig:continue_betamin1}
\end{figure}

\subsubsection{Limiting straight-line orbits}
\label{sec:line-orbit}

Having identified the scaling behavior of the periodic orbit approaching the locus \textbf{M}, we now examine candidate solutions for the respective limits. Since both amplitude and period grow beyond bound as either $\alpha$ increases to $0$ for fixed $\beta < -1$, or  $\beta$ increases to $-1$ for fixed $0 \leq \alpha < 1$, the length of the orbit segment where the derivative is virtually indistinguishable from $1-\alpha$ grows beyond bound as well. This suggests that the periodic behavior ceases in the limit: at the locus \textbf{M} the dynamics is an unbounded, linearly increasing segment $u_\ast$ connecting $-1$ and $\infty$, which we refer to as the limiting straight-line orbit.

It is advantageous to consider straight-line orbits for any $b \geq 0$. Substituting the ansatz $u_\ast(t)=k_{1}t+k_{2}$ into Eq.~\eqref{eq:sdDDE} implies
\begin{equation}
k_{1}(1+\beta(1- k_{1}b))=(k_{1}t+k_{2})(\alpha+\beta(1- k_{1})).
\label{eq:y-consistency-1}
\end{equation}
Any nontrivial solution of Eq.~\eqref{eq:y-consistency-1} requires
\[
\beta=\beta_\ast=-\frac{1-\alpha b}{1-b},\mbox{ and }k_{1}=k_\ast=1+\frac{\alpha}{\beta}=\frac{1-\alpha}{1-\alpha b}.
\]
Therefore, the straight-line orbit $u_\ast(t)=k_\ast t + k_2$ of slope $k_\ast $ solves Eq.~\eqref{eq:sdDDE} if and only if $\beta=\beta_\ast$; here the offset parameter $k_2$ is arbitrary and can be set to zero without loss of generality. In order to ensure that along $u_\ast$ the delay $\tau_\ast=1- u_\ast(\cdot-b)$ is well behaved, we require $\tau_\ast\geq 0$ (positivity) and $\frac{d}{dt}\tau_\ast\geq 0$ (monotonicity) for all $t$, which means equivalently $0\leq \alpha \leq \min\{1,1/b\}$.

In order to obtain a similar result for the specific limit $\alpha\to 0$ for fixed $\beta\leq -1$, we define $s=\alpha t$, $v(s)=\alpha u(t)$ and recast Eq.~\eqref{eq:sdDDE} with respect to the rescaled variables as
\begin{equation}
\label{eq:sdDDE-rescaled}
\alpha \frac{d}{ds} v(s)=\alpha v(s)+\beta v(s-\alpha-v(s-\alpha b)).
\end{equation}
Equation~\eqref{eq:sdDDE} is singularly perturbed. Setting $\alpha=0,$ we obtain $0=\beta v(s-v(s))$, which is trivially satisfied by the straight-line orbit $v_\ast(t)=t$ of slope $1$.

These two cases can be combined by defining the locus
\begin{equation}
\label{eq:locusL}
\mbox{\textbf{L}}=\left\{\left.(\alpha,\beta)\in\R^2 \ \right| \ \alpha\left(\beta+\frac{1-\alpha b}{1-b}\right)=0 \mbox{ with } 0\leq \alpha \leq 1,\,\beta\leq -\frac{1}{1-b} \right\},
\end{equation}
of straight-line orbits in $(\alpha,\beta)$-plane for any $b$. More specifically, for  $0\leq b <1$ the locus \textbf{L} consists of two parts: the line segment connecting the point \textbf{DZ} at $(\alpha,\beta) = (1,-1)$ with the point $(\alpha,\beta) = (0,-\frac{1}{1-b})$, which is a corner where \textbf{L} continues as the half-line given by $\alpha=0$ and $\beta \leq -1/(1-b)$. In particular, the loci \textbf{L} and \textbf{M} coincide for $b=0$, which is consistent with our observation that then the respective straight-line orbit reflects (at least part of) the limit object as the periodic orbit approaches the locus \textbf{M} = \textbf{L}.

\section{Bifurcation analysis of DDE with delayed state dependence}
\label{sec:bif-analysis}

We now study the effects of a delay $b>0$ in the state-dependent delay term of Eq.~\eqref{eq:sdDDE}. We do so by incrementally increasing the parameter $b$ from zero, thus, slowly `switching on' the delay in the state-dependency $\tau=1+u(\cdot-b)$, to determine how the bifurcation diagram in the $(\alpha,\beta)$-plane changes in the process. Nonzero $b$ introduces additional complications for defining solutions, which we evade by requiring positivity of the combined delay term $\tau$, or equivalently $u\geq-1$ uniformly, which is captured by the locus \textbf{M}. This restriction is very natural in the context of control theory because it ensures causality and, therefore, the physical relevance of solutions.

We first consider in Section~\ref{sec:smallb} the $(\alpha,\beta)$-plane for small values of $b$ up to $b=0.2$. This shows how the bifurcation diagram for $b=0$ changes with $b>0$; in particular, the loci \textbf{M} and \textbf{L} no longer coincide. We then consider in Section~\ref{sec:largerb} larger values of $b$ up to $b=0.5$ and show that this entails a qualitative change of the codimension-two point \textbf{DZ} and, hence, of the entire bifurcation diagram in the $(\alpha,\beta)$-plane. The associated degeneracy of the point \textbf{DZ} is studied by means of a center manifold analysis in Section~\ref{sec:degDZ}.

\subsection{Bifurcation diagram for small $\boldsymbol{b>0}$}
\label{sec:smallb}

\begin{figure}[t!]
	\centering
	\includegraphics[width=0.76\textwidth]{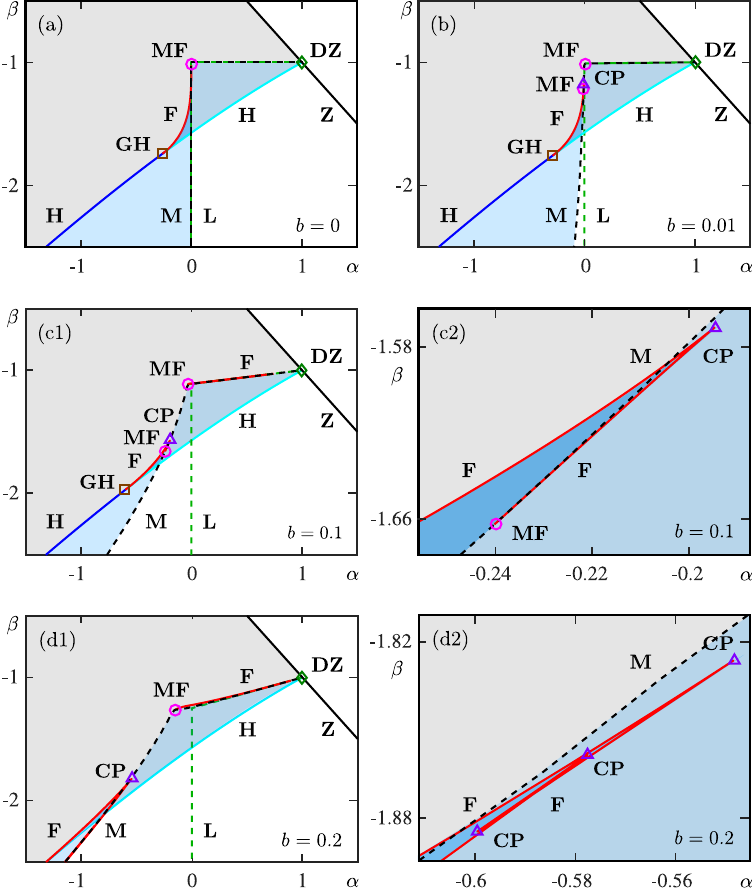}
	\caption{Bifurcation diagram of system~\eqref{eq:sdDDE} in the $(\alpha,\beta)$-plane for (a) $b=0$, (b) $b=0.01$, (c1) $b=0.1$  with enlargement (c2); and (d1) $b=0.2$  with enlargement (d2). Shown are the curves \textbf{Z} (black) of simple zero eigenvalue and \textbf{H} of Hopf bifurcations (blue when supercritical, cyan when subcritical) of $u=0$, the curve \textbf{F} (red) of fold bifurcation of periodic orbits, and the loci \textbf{M} (black dashed) and \textbf{L} (green dashed); also shown are the points \textbf{DZ} (green diamond) of double zero eigenvalues, \textbf{GH} (brown square) of generalized Hopf bifurcation, \textbf{CP} (purple triangle) of cusp bifurcation on \textbf{F}, and \textbf{MF} (magenta circle) of fold periodic orbit with minimum value $-1$. The stability region of $u=0$ is shaded gray, and blue shading indicates (co)existence of periodic orbits.}
\label{fig:Bifurcation-diagram}
\end{figure}

Figure~\ref{fig:Bifurcation-diagram} shows how the bifurcation diagram of system~\eqref{eq:sdDDE} in the relevant part of the $(\alpha,\beta)$-plane near the point \textbf{DZ} changes when $b$ is increased up to $b=0.2$. As the starting point, panel~(a) shows the bifurcation diagram for $b=0$ from Fig.~\ref{fig:stabdiag_b0}(b) over a slightly extended range; here the locus \textbf{L} given by Eq.~\eqref{eq:locusL} has also been added and labeled, showing that indeed \textbf{L} agrees with \textbf{M} in this case. When $b$ is increased, the curves \textbf{Z}, \textbf{H} and the point \textbf{DZ} are unaffected as they are determined by the linear stability analysis alone. However, on the nonlinear level $b$ does have an effect and the bifurcation diagram changes as soon as $b$ is increased from $0$. As Fig.~\ref{fig:Bifurcation-diagram}(b) for $b=0.01$ shows, there are two qualitative changes. Firstly, the locus \textbf{M}  that bounds the region of periodic orbits (with positive state-dependent delay) has deviated from the locus \textbf{L} along which one finds straight-line orbits; this is visible chiefly along their vertical sections, and \textbf{M} and \textbf{L} are still very close together along their horizontal sections. While \textbf{L} is plotted readily from Eq.~\eqref{eq:locusL}, the locus \textbf{M} can be continued reliably as a periodic orbit in two parameters subject to the additional condition that the minimum value along its profile equals $-1$. Secondly, the corner point \textbf{MF} of \textbf{L} = \textbf{M} in panel~(a) for $b=0$ has split up in panel~(b) for $b=0.01$ into two points labeled \textbf{MF}, with the lower of the two now being the end point of the curve \textbf{F} of fold bifurcations of periodic orbits that emerges from the generalized Hopf bifurcation point \textbf{GH}; moreover, there is now a cusp point \textbf{CP} on the curve \textbf{F}, close to where it reaches \textbf{MF}.

When $b$ is increased further, these qualitative changes become more visible, as is shown in Fig.~\ref{fig:Bifurcation-diagram}(c1) for $b=0.1$. The deviation of the locus \textbf{M} from the vertical part of \textbf{L} is now very pronounced, leading to considerable shrinking of the region where one finds stable periodic orbits (that emerge from the supercritical part of \textbf{H}). Since the points \textbf{GH}, \textbf{CP} and \textbf{MF} have also moved to lower values of $\beta$, the region with two periodic orbits has become vanishingly small. The enlargement panel~(c2) shows that the curve \textbf{F}, when continued from \textbf{GH}, crosses \textbf{M} and has a cusp point \textbf{CP} before the curve \textbf{F} ends on \textbf{M} at the point \textbf{MF}. Notice in Fig.~\ref{fig:Bifurcation-diagram}(c1) that \textbf{M} is still very close to the still almost horizontal slanted segment of \textbf{L}. Nevertheless, close to this segment of \textbf{M} we were able to detect fold bifurcations of periodic orbits that have both large period and large amplitude. The continuation of such folds is very challenging numerically, but we were able to find the fold curve \textbf{F} shown as a new feature of the bifurcation diagram for $b>0$. The curve \textbf{F} emerges from the upper point \textbf{MF} on the locus \textbf{M} and tracks the almost horizontal segment of \textbf{M} very closely. While the continuation stops somewhat short of this point when $b=0.1$, we suspect that \textbf{F} extends all the way to the point \textbf{DZ}. This is confirmed in Fig.~\ref{fig:Bifurcation-diagram}(d1) for the larger value $b=0.2$, for which we were able to continue the curve \textbf{F} all the way from \textbf{MF} to \textbf{DZ}. In spite of the fact that the respective fold bifurcations of periodic orbits are practically impossible to identify and continue when $b$ is very small, due to their extremely large periods, we therefore conclude that this additional curve \textbf{F} exists for any $b>0$. This numerical observation is confirmed in Sec.~\ref{sec:degDZ} via a normal form analysis of the point \textbf{DZ}. Notice further in Fig.~\ref{fig:Bifurcation-diagram}(d1) that the points \textbf{GH} on the curve \textbf{H} and \textbf{MF} on the curve \textbf{M} have moved to small values of $\beta$ outside the region of the $(\alpha,\beta)$-plane shown here. Otherwise, the bifurcation diagram for $b=0.2$ appears to be qualitatively the same as that for $b=0.1$. However, this is not quite the case: as can be seen in the enlargement panel~(d2), close to the cusp bifurcation point \textbf{CP}, there are now two additional cusp points, also labeled \textbf{CP}, which are the result of a transition through a swallowtail bifurcation (degenerate cusp bifurcation) \cite{Arnold92,Guc-Hol-90,Kuz-04-Book}. Notice this extra feature is on a very small scale in the $(\alpha,\beta)$-plane and hardly changes the boundary of the region with two periodic orbits.

An important feature of Fig.~\ref{fig:Bifurcation-diagram} is the deviation of the locus \textbf{M} from the locus \textbf{L}, and we proceed by discussing the behavior of periodic orbits as \textbf{M} is approached for the case that $b=0.1$. Here, similar to our study for $b=0$ in Sec.~\ref{sec:saw-tooth-limits},  we first consider the more vertical section of \textbf{M}, which is quite far from \textbf{L}, and then the almost horizontal section of \textbf{M}, which is still very close to \textbf{L}.

\subsubsection{Periodic orbits when approaching the more vertical part of  {\normalfont\textbf{M}}}
\label{sec:verticalM}

\begin{figure}[t!]
	\centering
    \includegraphics[width=0.65\textwidth]{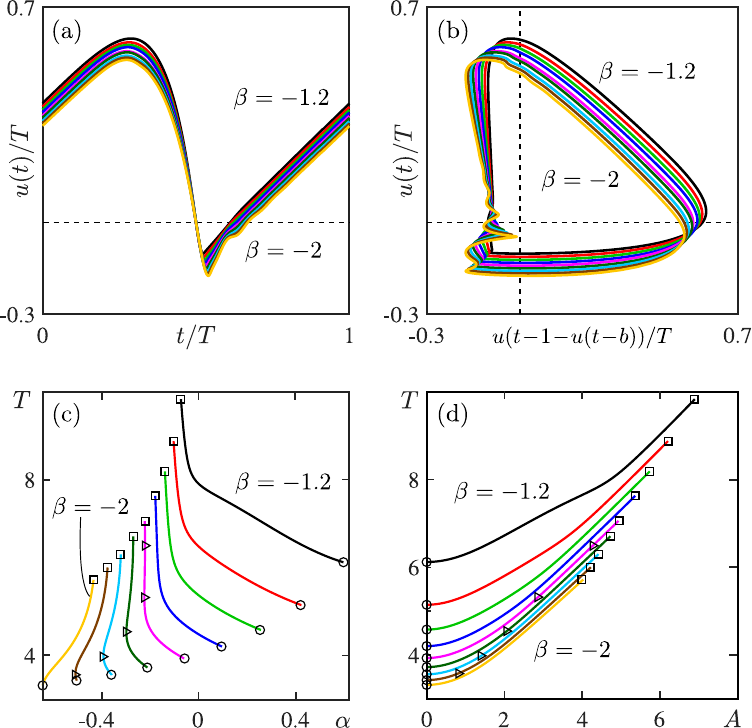}\hspace{0.35cm}
	\caption{Periodic orbits of the system~\eqref{eq:sdDDE} for $b=0.1$ as $\alpha$ approaches the locus \textbf{M} for the different fixed values of $\beta=-1.2$ down to $\beta=-2$. Panel~(a) shows their (rescaled) profiles at \textbf{M} synchronized to their fast segments, and panel~(b) shows them in the $(u(t-1-u(t-b))/T,u(t)/T)$-plane; compare with Fig.~\ref{fig:limit_alpha0}. Panel~(c) shows the continuations in $\alpha$ of the periodic orbits, represented by the period $T$, from the Hopf bifurcation \textbf{H}  (circles), possibly via fold points (triangles) to the locus \textbf{M} (squares); the same branches are shown in panel~(d) as a function of the amplitude $A$ of the periodic orbits; compare with Fig.~\ref{fig:continue_alpha0}.}
\label{fig:scaling-b01-constbeta}
\end{figure}

Figure~\ref{fig:scaling-b01-constbeta} shows the fate of nine distinct stable periodic orbits of system~\eqref{eq:sdDDE} for $b=0.1$ as $\alpha$ approaches the locus \textbf{M} for different fixed values of $\beta\leq-1.2$. Panels~(a) and~(b) show the periodic orbits at the locus \textbf{M} as profiles and in projection onto the $(u(t),u(t-1-u(t)))$-plane, respectively, where all axes have again been rescaled by $T$. Comparison with Fig.~\ref{fig:limit_alpha0} for $b=0$ shows that these orbits no longer have a very distinctive saw-tooth shape: while they all feature an almost linearly increasing segment and a sharp corner when reaching their respective minimum value, the periodic orbits are still quite far from the limiting straight-line orbit and, hence, still clearly distinguishable. Moreover, we observe in Fig.~\ref{fig:scaling-b01-constbeta}(a) and~(b) that decreasing $\beta$ leads to secondary oscillations along the periodic orbits after the switch to negative values; these oscillations are more pronounced the smaller $\beta$. One-parameter continuation of these periodic orbits in $\alpha$, from where they are created at the Hopf bifurcation curve \textbf{H} all the way to the locus \textbf{M}, are shown in panels~(c) and~(d). As was the case for $b=0$ in Fig.~\ref{fig:continue_alpha0},  there are fold points on some of the branches of periodic orbits in Fig.~\ref{fig:scaling-b01-constbeta}(c), namely for $\beta=-1.6$ down to $\beta=-1.9$. More specifically, for $\beta$ above the cusp point \textbf{CP} on the curve \textbf{F} at $\beta\approx-1.57$ there or no folds; in between \textbf{CP} and the point \textbf{MF} on the curve \textbf{M} at $\beta\approx-1.66$ there are two fold points; in between \textbf{MF} and the point \textbf{GH} on the Hopf bifurcation curve \textbf{H} at $\beta\approx-1.93$ there is one fold point; and for $\beta$ below \textbf{GH} there are again no folds on the branch of periodic orbits; compare with Fig.~\ref{fig:Bifurcation-diagram}(c1) and~(c2). Notice further that, as a result, the limiting periodic orbit on \textbf{M} in Fig.~\ref{fig:scaling-b01-constbeta} is stable below \textbf{MF} and unstable above \textbf{MF}. In either case, the period and the amplitude of the periodic orbit reaches a finite limit, that is, they are bounded on \textbf{M}. This is owing to the fact that \textbf{M} and \textbf{L} no longer coincide and illustrates comprehensively that the straight-line orbit is no longer a good approximation for the orbits along this part of the locus \textbf{M} (with $\beta< -1.1$).

\subsubsection{Periodic orbits when approaching the almost horizontal part of  {\normalfont\textbf{M}}}
\label{sec:horizontalM}

\begin{figure}[t!]
	\centering
    \includegraphics[width=0.65\textwidth]{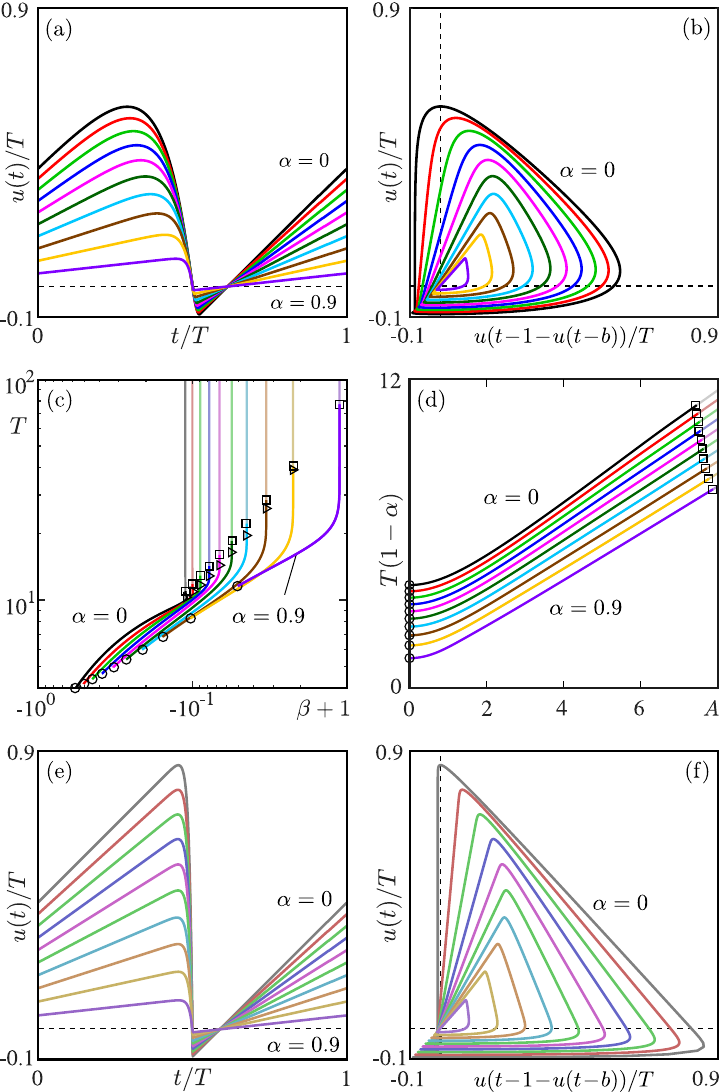}\hspace{0.35cm}
	\caption{Periodic orbits of system~\eqref{eq:sdDDE} for $b=0.1$ as $\beta$ approaches \textbf{M} for ten values of $\alpha=0$ up to $\alpha=0.9$. Panels~(a) and~(b) show the periodic orbits at \textbf{M} as (rescaled) profiles and in the $(u(t-1-u(t-b))/T,u(t)/T)$-plane, respectively; compare with Fig.~\ref{fig:limit_betamin1}. Panels~(c) and~(d) show their period $T$ and amplitude $A$ when continued in $\beta$ from \textbf{H}  (circles), via fold points (triangles) to \textbf{M} (squares) and beyond (lighter curves) up to $T=100$. Panels~(e) and~(f) show the periodic orbits with $T=100$ as profiles and in the $(u(t-1-u(t-b))/T,u(t)/T)$-plane, respectively.}
\label{fig:scaling-b01-constalpha}
\end{figure}

Figure~\ref{fig:scaling-b01-constalpha} shows ten distinct unstable periodic orbits of system~\eqref{eq:sdDDE} for $b=0.1$ as $\beta$ approaches the locus \textbf{M} for different fixed values of $\alpha \in [0,0.9]$. Panels~(a) and~(b) again show these periodic orbits at \textbf{M} (where they have a minimum of $-1$) as profiles and in the $(u(t-1-u(t))/T,u(t)/T)$-plane, respectively. The periodic orbits all feature an almost linearly increasing segment whose slope scales like $1-\alpha$, as well as sharp corners at their minima ---  yet are clearly still quite far from being of saw-tooth shape; compare with Fig.~\ref{fig:limit_betamin1}. This is in spite of the fact that the respective parts of the loci \textbf{M} and \textbf{L} are very close to each other in the $(\alpha,\beta)$-plane of Fig.~\ref{fig:Bifurcation-diagram}(c1). Figure~\ref{fig:scaling-b01-constbeta}(c) shows the associated one-parameter continuations of these periodic orbits in $\beta$, from where they are created at the Hopf bifurcation curve \textbf{H} to the locus \textbf{M}. Notice from Fig.~\ref{fig:Bifurcation-diagram}(c1) that there is a fold curve \textbf{F} above \textbf{M}, very close to its almost horizontal part; see also Fig.~\ref{fig:Bifurcation-diagram}(d1). Hence, the unstable periodic orbits, which bifurcate from the subcritical part of \textbf{H}, become stable just before \textbf{M} is reached. Indeed, we found fold points very close to \textbf{M} on the branches for $\alpha = 0$ up to $\alpha = 0.8$ in Fig.~\ref{fig:Bifurcation-diagram}(c1); for $\alpha = 0.9$, on the other hand, we do not show a fold point very close to \textbf{M} because its detection is no longer reliable since the branch is practically vertical. Panel~(d), showing the same branches as a function of the amplitude $A$, provides evidence that the expected scaling behavior $A \sim T(1-\alpha)$ between amplitude and period still applies, although both $A$ and $T$ only reach finite values at \textbf{M}. 

Each of the branches of (now stable) periodic orbits can actually be continued in $\beta$ past the locus \textbf{M}, in spite of the fact that the periodic orbits now contain a segment where the delay is positive; this is possible because periodic orbits are continued as solutions of a boundary value problem. The lighter curves in Fig.~\ref{fig:scaling-b01-constbeta}(c) and~(d) are the continuations of the ten branches computed up to a period of $T=100$. Notice that these extended branches are practically vertical in panel~(c), meaning that the parameter $\beta$ hardly changes, which is due to the very small distance between the loci \textbf{M} and \textbf{L}. The periodic orbits with $T=100$ are shown panels~(e) and~(f) as time series and in the $(u(t-1-u(t))/T,u(t)/T)$-plane, respectively; this shows that they are much more saw-tooth shaped and, hence, closer to the locus \textbf{L}; compare with Fig.~\ref{fig:limit_betamin1}. We note here that these solutions past \textbf{M} are not directly physically relevant, but nevertheless serve to build intuition regarding the limiting shape of periodic orbits close to \textbf{L}.

\subsection{Qualitative changes near \textbf{DZ} for larger values of $\boldsymbol{b}$}
\label{sec:largerb}

\begin{figure}[t!]
	\begin{raggedright}
		\centering
		\includegraphics[width=0.76\textwidth]{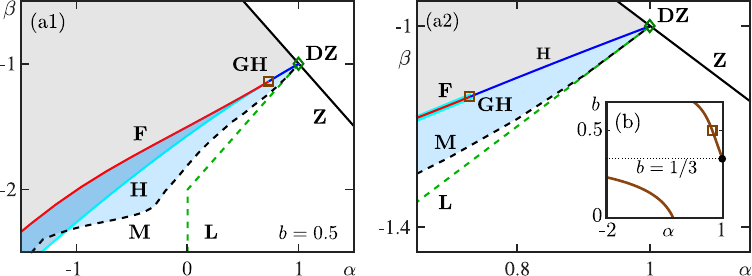}
		\par\end{raggedright}
	\caption{Bifurcation diagram of system~\eqref{eq:sdDDE} in the $(\alpha,\beta)$-plane (a1) for $b=0.5$, and (a2) an enlargement near the point \textbf{DZ}; compare with Fig.~\ref{fig:Bifurcation-diagram}. The inset panel~(b) shows the branches of generalized Hopf bifurcation points \textbf{GH} along the Hopf curve \textbf{H} in the $(\alpha,b)$-plane; the brown square highlights the point \textbf{GH} for $b=0.5$ shown in panels~(a1) and~(a2).}
\label{fig:Bifurcation-diagram-extended}
\end{figure}

We have observed in the previous section that the curve \textbf{F} and the loci \textbf{M} and \textbf{L} are very close together near the point \textbf{DZ}. Moreover, we know from Eq.~\eqref{eq:locusL} that the slope $b/(1-b)$ of \textbf{L} near \textbf{DZ} increases with increasing $b$, while the curve \textbf{H} does not change. Therefore, the curves \textbf{L} and \textbf{H} exchange their relative positions near \textbf{DZ} as $b$ is decreased further. We now discuss and illustrate in Fig.~\ref{fig:Bifurcation-diagram-extended} what this change implies for the curves \textbf{F} and \textbf{M} and the bifurcation diagram in the $(\alpha,\beta)$-plane more generally.

Panel~(a1) of Fig.~\ref{fig:Bifurcation-diagram-extended} shows the bifurcation diagram for $b=0.5$, and panel~(a2) is an enlargement near the point \textbf{DZ}. There clearly has been a qualitative change of the bifurcation diagram. Both the loci \textbf{L} and \textbf{M} now lie to the lower side of the curve \textbf{H}. Moreover, there is now a second point \textbf{GH} on \textbf{H}, where the criticality of the Hopf bifurcation changes. Note that \textbf{GH} is now the end point of the fold bifurcation curve \textbf{F}, rather than the point \textbf{DZ}; compare with Fig.~\ref{fig:Bifurcation-diagram}(d1). Between \textbf{GH} and \textbf{DZ} the bifurcating periodic orbit is stable and disappears (for decreasing $\beta$) at the locus \textbf{M}; see Fig.~\ref{fig:Bifurcation-diagram-extended}(a2). Notice the considerable region with two periodic orbits in panel~(a1), bounded by the subcritical part of \textbf{H} and the fold bifurcation curve \textbf{F}. The inset panel~(b) shows the $\alpha$-values of both generalized Hopf bifurcation points \textbf{GH} against $b$. The branch that emerges from $b=0$ is the point \textbf{GH} we discussed previously; it moves towards larger and larger negative values for $\alpha$ along the Hopf bifurcation curve \textbf{H} as $b$ increases, and it is outside the shown parameter range already for $b=0.2$. However, there is a second branch, which emerges from $\alpha=1$, that is, from the point \textbf{DZ}, at $b=1/3$; its instance for $b=0.5$ is highlighted by a square.

\section{Behavior near the point \textbf{DZ}}
\label{sec:degDZ}

Our numerical investigation clearly shows that there is a change in the nature of the degenerate transcritical bifurcation at \textbf{DZ} when $b=1/3$, and that this is the reason for the qualitative change of the bifurcation diagram illustrated in Fig.~\ref{fig:Bifurcation-diagram-extended}(a1) and~(a2). This motivates a center manifold analysis of the codimension-two point \textbf{DZ} here. The double-zero singularity at the linear level has already been identified in \cite{walther2006bifurcation}, where its consequences are studied for the case of a DDE with constant delay and a non-smooth instantaneous term of the form
  \begin{align}\label{walther:dde}
    a^{-1} \dot u(t)&=(1-|u(t)|)u(t)-u(t-1)\mbox{.}
  \end{align}
The analysis of this system in \cite{walther2006bifurcation} established that a family of periodic orbits emerges from $u=0$ at \textbf{DZ}, with period going to infinity while the amplitude decreases to zero; only the parameter $a$ was varied, and the non-smoothness of \eqref{walther:dde} precluded expansion and truncation. In our case, the semiflow also lacks smoothness due to state dependence, but we may still apply expansion techniques near \textbf{DZ}. Furthermore, our nonlinearity is infinitely degenerate, leading to a line of equilibria in \textbf{DZ}. Our goal is to identify and unfold the codimension-three point of the DDE \eqref{eq:sdDDE} at $(\alpha,\beta,b)=(1,-1,1/3)$, which we refer to as \textbf{DZGH}. To this end, we write system~\eqref{eq:sdDDE} as 
\begin{align*}
f: C & \to \R \\
u & \mapsto \alpha u(0)+\beta u(-1-u(-b))
\end{align*}
with the arbitrarily often mildly differentiable right-hand side given by $f$; see \cite{Hartung2006,Sieber12,Sieber2017} for the definition of mild differentiability. We use the convention that $C^k = C^k([-\tmax,0;\R)$ is the space of $k$ times continuously differentiable functions on $[-\tmax,0]$. The equilibrium $u=0$ of the \eqref{eq:sdDDE} is degenerate for $\alpha=1$ and $\beta=-1$.  More precisely, along \textbf{Z} we have $\alpha=-\beta$ and Eq.~\eqref{eq:sdDDE} has a one-parameter family of equilibria, $x(t)=x_\mathrm{eq}$ with arbitrary $x_\mathrm{eq}$. Consequently, the linearization 
\begin{align}\label{lindde}
x'(t) =\alpha Lx_t \quad \mbox{with} \quad L:C & \to \R \\
x & \mapsto x(0)-x(-1)\mbox{} \nonumber
\end{align}
of Eq.~\eqref{eq:sdDDE} about $x=0$ (for arbitrary $\alpha$ and $\beta=-\alpha$) has a first trivial eigenvalue $0$ with eigenfunction $x_1(t)=1$. In what follows, we restrict attention to the degenerate case $\alpha=1$, where the characteristic function  $\chi(\lambda)=\lambda-(1-\exp(-\lambda))$ in $x=0$ has a double root $0$; in other words, the linear operator associated with the linear DDE \eqref{lindde},
\begin{align}
A: C^1 & \to C^0 \nonumber \\
x & \mapsto x' \quad \mbox{with domain} \ \ 
D(A) =\{x\in C^1:x'(0)=Lx\}\mbox{,}
\label{eq:derA}
\end{align}
has an eigenvalue $0$ with algebraic multiplicity $2$ at the point \textbf{DZ}. The curve \textbf{H} of Hopf bifurcations given in Eq.~\eqref{eq:C0} emerges from $(\alpha,\beta)=(1,-1)$ and has the expansion 
\begin{align*}
\alpha=1-\frac{\omega^2}{3} + O(\omega^4) \quad \mbox{and} \quad \beta =-1-\frac{\omega^2}{6} + O(\omega^4)\mbox{}
\end{align*}
for small frequency $\omega$. 

We now shift \textbf{DZ} to the origin of the parameter plane and align the loci \textbf{H} and \textbf{Z} locally with the coordinate axes. This is achieved by introducing parameters $(p,q)$ that are related to $(\alpha,\beta)$ via
\begin{align}
\label{eq:alphabeta_to_pg}
	\begin{bmatrix}
		\alpha\\\beta
	\end{bmatrix}=
	\begin{bmatrix}
		\phantom{-}1\\-1
	\end{bmatrix}+
	\begin{bmatrix}
		\phantom{-}1/3&1/2\\-1/6&1/2
	\end{bmatrix}
	\begin{bmatrix}
		p\\ q
	\end{bmatrix}\mbox{,}
\end{align}
such that $(p,q)$ are small near the degenerate point \textbf{DZ}. The coordinates $x_c$ on the center manifold are in $\R^2$, and map into the space $C^1$ via the graph $h(x_x;p,q)=Bx_c+O(2)$, where $O(2)$ denotes second order terms in $x_c$ and $(p,q)$. Thus, we expect that the dynamics near \textbf{DZ} is governed by a two-dimensional ODE on this center manifold. The linear transformation
\begin{equation}
\label{eq:xctrafo}
x_c=
\begin{bmatrix}
1&1/3\\0&1
\end{bmatrix}
y_c
\end{equation} 
%
transforms the ODE on the center manifold into a second-order differential equation such that $\dot y_{c,1}=y_{c,2}$ for expansions at all orders $j$. As a result we obtain the following.

\begin{lemma}
\label{thm:cmf2}
To second order, the resulting ODE for $y=y_{c,1}$ on the center manifold of the point \textbf{DZ} is
\begin{equation}
\label{eq:cmf2}
\ddot y= py+q\dot y+2y\dot y+(2/3-2b)\dot y^2+O(3)\mbox{.}
\end{equation}
\end{lemma}
The derivation and terms up to order five are provided in Appendix~\ref{app:ex} and the supplementary material with this paper. In particular, for expansions up to all orders the right-hand side is zero when $p=\dot y=0$, meaning that the line of equilibria (the $y$-axis where $\dot{y}=0$) at $p=0$ is always present (it is stable for $y<-q$ and unstable for $y>-q$).

\subsection{Qualitative analysis of the truncated ODE on the center manifold}
\label{sec:cmfdynamics}

The nonlinear part of the right-hand side of \eqref{eq:sdDDE} takes the form
\begin{align*}
f_{(p,q)}(u)=(p/6 - q/2)u(0) + (p/6 + q)u(-1-u(-b))
\end{align*}
in the transformed parameters $p$ and $q$. The two-dimensional local center manifold of \eqref{eq:sdDDE} near \textbf{DZ} of the semiflow on the manifold
\begin{align*}
\{u\in C^1:u'(0)=u(0)-u(-1)+f_{(p,q)}(u)\},
\end{align*}
which is known to exist for $p\approx0$ ($\alpha\approx1$),  $p\approx0$ ($\beta\approx-1$) and $u\approx0$, is a graph over the basis functions $B_1(\theta)=1$ and $B_2(\theta)=\theta$ \cite{Stu-12}.  As \cite{Hum-Dem-Mag-Uph-12,Sieber2017} laid out, all formal expansion coefficients of this center manifold $h:\R^2\times\R^2\mapsto C$ are well defined to arbitrary order. Appendix~\ref{app:ex} gives details of the resulting linear systems of algebraic equations; we also show there that higher-order derivatives of the mildly differentiable right-hand side $f_{(p,q)}$ are always applied to $C^\infty$ history segments (which are prior expansion terms of the center manifold and, therefore, solutions to linear ODEs).

\begin{figure}[t!]
  \centering
  \includegraphics[width=0.85\textwidth]{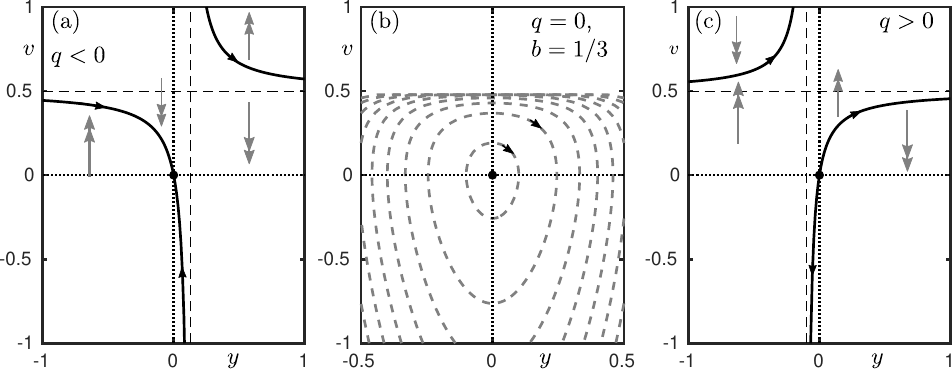}
	\caption{Illustration of the truncated dynamics on the center manifold near \textbf{DZ} when the parameter $p$ is small and negative ($-1\ll p<0$). Panels (a) and (c) show the slow manifold $v=y/(q+2y)$, which is a small perturbation of the line of equilibria at $p=0$ and the fixed point at $(0,0)$. Panel (b) shows the level curves of the potential $V$ in \eqref{eq:consquant} for the case $q=0$ and $b=1/3$, where the truncated system \eqref{eq:cons2} is conservative.}
\label{fig:slowmanifold}
\end{figure}

\smallskip
\noindent
\emph{Dynamics near a slow manifold for $-1\ll p<0$.} For $p=0$, one may apply the analysis by Liebscher in~\cite{liebscher15}. A family of arcs in the half-plane $\{\dot y<0\}$ connects each point $(y_1,0)$ on the part $\{y>-q, \dot y=0\}$ of the line of equilibria with one point $(y_2,0)$ on $\{y<-q, \dot y=0\}$. Since for $p=0$ the system is degenerate (having a line of equilibria), it is singularly perturbed for $p\neq0$ (in particular, for $p<0$). Speeding up time by a factor of $-p$ clarifies the slow-fast structure in the second-order truncation
\begin{align}\label{eq:cmf_rescaled}
	y' = v\mbox{,} \quad (-p)v' = -y+qv+2yv-(2/3-2b)pv^2\mbox{,}
\end{align}
where $v=y'$ and $v'$ denote the derivatives of $y$ with respect to the new time. 

In the $(y,v)$ coordinates of the center manifold the slow manifold is $v=y/(q+2y)$, and this is shown in Fig.~\ref{fig:slowmanifold}. The center manifold is transversally stable for $y<-q/2$ and unstable for $y>-q/2$. The line $y=-q/2$ is a first-order pole of the manifold for $q\neq0$ and a singular point for $q=0$. For $q=0$ the slow manifold has a singular branch $y=0$, intersecting the regular slow manifold $v=1/2$ at $(y,v)=(0,1/2)$. Moreover, for $b=1/3$ the truncated system \eqref{eq:cons2} is conservative and the equilibrium $(0,0)$ is a center; see Fig.~\ref{fig:slowmanifold}(b). When $q \neq 0$ then the equilibrium $(0,0)$ lies on the stable slow manifold for $q<0$ and on the unstable slow manifold for $q>0$; see panels~(a) and~(c). The dynamics along the slow manifold, given by $y'=y/(q+2y)$, leads to monotone increasing $y$ outside of the intervals $[0,-q/2]$ for $q<0$ and $[-q/2,0]$ for $q>0$, respectively, (where $y$ is decreasing).

\smallskip
\noindent
\emph{Conserved quantity for $q=0$ and $b=1/3$.} The Hopf bifurcation is located at $q=0$ and $p<0$. The second-order truncated ODE \eqref{eq:cmf2} has a degeneracy if $b=1/3$ and $q=0$. For these parameters, the second-order truncated rescaled ODE \eqref{eq:cmf_rescaled} reads
\begin{align}\label{eq:cons2}
	y' =v\mbox{,} \quad  (-p)v' = -y+2yv\mbox{}
\end{align}
and has the conserved quantity 
\begin{align}\label{eq:consquant}
V(y,v)=-\log(pv-p/2)-2v - 2y^2/p,
\end{align}
which has a center at $(y;v)=(0;0)$ for $p<0$. Thus, the Hopf bifurcation is subcritical for $b<1/3$, supercritical for $b>1/3$ and infinitely degenerate in the truncated system for $b=1/3$ and $q=0$. Figure~\ref{fig:slowmanifold}(b) shows the contours of the quantity $V$.

\subsection{Numerical bifurcation diagrams of center manifold expansion up to order five}
\label{sec:cmfbifs}

\begin{figure}[t!]
	\centering
	\includegraphics[width=0.76\textwidth]{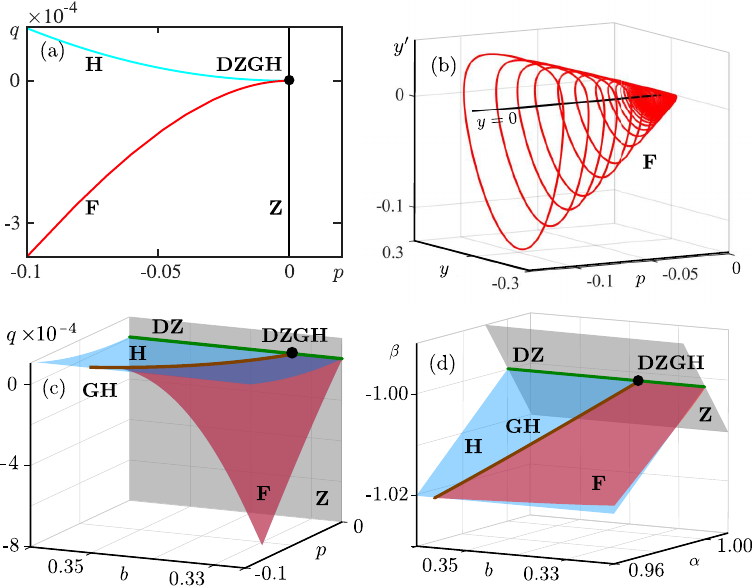} 	\caption{Bifurcation diagrams near the point \textbf{DZGH}. Panel~(a) is the two-parameter bifurcation diagram in the $(p,q)$-plane of the fifth-order expansion on the center manifold for $b=1/3$; shown are the curves \textbf{Z} of zero-eigenvalue equilibria, \textbf{H} of Hopf bifurcations and \textbf{F} of fold periodic orbits, which intersect at the point \textbf{DZGH}. Panel~(b) shows the family of fold periodic orbits along \textbf{F} in the $(y, y')$-plane as a function of the parameter $p$ (with $q$ varying accordingly, but not shown). The three-parameter bifurcation diagram near \textbf{DZGH} of the fifth-order expansion in $(b,p,q)$-space near $(1/3,0,0)$ is shown in panel~(c), and that of the full DDE \eqref{eq:sdDDE} in $(b,\alpha,\beta)$-space near $(1/3,1,-1)$ is shown in panel~(d); both consist of the surfaces \textbf{Z}, \textbf{H} and \textbf{F}, and the curves \textbf{DZ} and \textbf{GH} that meet at the point \textbf{DZGH}.}
\label{fig:bif3d}
\end{figure}

The second-order expansion \eqref{eq:cmf2} does not unfold the degeneracy of the Hopf bifurcation at $b=1/3$, indicated by the conserved quantity \eqref{eq:consquant}, but higher order terms in the center manifold expansion \eqref{eq:cmf2} do. We find that the fifth-order expansion results in a surface \textbf{F} of folds (saddle-node bifurcations) of periodic orbits that appear from a curve \textbf{GH} of generalized Hopf bifurcations. This is illustrated in Fig.~\ref{fig:bif3d} in different ways. Panel~(a) shows the 
two-parameter bifurcation diagram in the $(p,q)$-plane (of the rescaled parameters) for $b=1/3$ of the fifth-order expansion of the ODE on the center manifold; here the curves \textbf{H} of Hopf bifurcations and \textbf{F} of fold periodic orbits meet at the point \textbf{DZGH} at $(p,q) = (0,0)$ on the curve \textbf{Z} of zero-eigenvalue equilibria. The fold periodic orbits along \textbf{F} (of the fifth-order expanded and truncated ODE on the $(y, y')$-plane) are shown in Fig.~\ref{fig:bif3d}(b) in $(p,y, y')$-space, which illustrates how they shrink to a point as $p \to 0$ and the point \textbf{DZGH} is approached; note that $q$ varies accordingly along the curve \textbf{F}, but is not shown in this illustration. 

The unfolding of the point \textbf{DZGH} by the fifth-order ODE on the center manifold in $(b,p,q)$-space is shown in Fig.~\ref{fig:bif3d}(c). It illustrates how the curve \textbf{GH} of generalized Hopf bifurcations emerges from \textbf{DZGH} at $(b,p,q) = (1/3,0,0)$; together with the part of the curve \textbf{Z} with $b<1/3$, the curve \textbf{GH} forms the boundary of the surface \textbf{F} of fold periodic orbits. Figure~\ref{fig:bif3d}(d) shows for comparison the same surfaces and curves of the full DDE \eqref{eq:sdDDE} in $(b,\alpha,\beta)$-space near the point \textbf{DZGH} at $(b,\alpha,\beta) = (1/3,1,-1)$. Clearly, the fifth-order ODE on the center manifold captures the local behavior near the point \textbf{DZGH} that is observed in the DDE. In particular, the asymptotics of the folds or periodic orbits on \textbf{Z} can be captured with a truncation of the ODE at a suitable order; hence, terms of the ODE on the center manifold that are beyond all orders, as analyzed by \cite{JMPRNP94}, do not play a role for the fold of periodic orbits.

\section{Conclusions and outlook}
\label{sec:conclusions}

We investigated the interplay of instantaneous and delayed feedback in a scalar, first-order DDE; here, the delayed feedback term depends linearly on the state, but the state itself may enter the delay term in a delayed way. The system is linear in the absence of state dependence, and we studied the effects of the nonlinearity induced by state depencence of the delay on the bifurcation diagram in the two-parameter plane of the instantaneous feedback strength $\alpha$ and delayed feedback strength $\beta$. The system exhibits periodic solutions that are generated via Hopf bifurcations, and we showed that they cease to exist (become unphysical) when the state-dependent delay becomes advanced, or equivalently, when the periodic orbit attains its minimum value at a certain threshold value (here rescaled to $-1$).

We showed that, when the state dependency is instantaneous, the latter occurs along a well-defined curve \textbf{M} in the $(\alpha,\beta)$-plane, approaching which the periodic orbits develop a distinct saw-tooth shape. This allowed us to present a complete bifurcation diagram, with curves \textbf{Z} of zero-eigenvalue equilibria, \textbf{H} of Hopf bifurcation, \textbf{F} of folds of periodic orbits and \textbf{M} of limiting saw-tooth periodic orbits bounding the existence and stability regions of the basic equilibrium and of the bifurcating periodic solutions. Further important features of this bifurcation diagram are a generalized Hopf bifurcation point \textbf{GH} where the Hopf bifurcation curve changes criticality and a point \textbf{DZ} where the basic equilibrium has a double eigenvalue zero.

The case of instantaneous state dependence formed the basis for our investigation of the dynamics that arise when the state-dependent feedback term is itself subject to a delay $b$. We found that instantaneous state dependence is special in the sense that `switching on' this additional delay by considering $b > 0$ immediately changes the bifurcation diagram in the $(\alpha,\beta)$-plane. To show how this happens, we employed advanced numerical continuation techniques for DDEs to find the respective two-parameter bifurcation diagram for fixed values of $b$. In particular, we found that the curve \textbf{M} where periodic orbits become unphysical changes rapidly with the additional delay $b$ and, moreover, gives rise to a further curve \textbf{F} of folds of periodic orbits, as well as a pair of points \textbf{MF} where the two curves \textbf{F} end on the curve \textbf{M}. Furthermore, the bifurcation diagram in the $(\alpha,\beta)$-plane changes qualitatively at $b = 1/3$, where the point \textbf{DZ} has an additional degeneracy that gives rise to second generalized Hopf point \textbf{GH} on the curve \textbf{H}. A center manifold analysis, by means of computing the expansion of the ODE on the center manifold up to order five, confirmed this result.

The results presented here for this prototypical example are expected to be of interest more widely, beyond their significance for the theory of state-dependent DDEs. They further elucidate the capacity of state dependency alone to give rise to nonlinear dynamics in any model under consideration, irrespective of other sources of nonlinearity. Furthermore, the results suggest that instantaneous state dependence may be special mathematically: as we have shown, the observed dynamics may change significantly when there is even a small delay in the state-dependent feedback term. Additional qualitative changes may occur when this additional delay becomes larger, as we showed with the transition through the degenerate point \textbf{DZGH} at $b = 1/3$; the study of further qualitative changes of the bifurcation diagram in the $(\alpha,\beta)$-plane with increasing $b$ is the subject of ongoing research.

An interesting direction for future research is the in-depth analysis of periodic orbits as the delay becomes advanced. We observe that, for $b>0$, secondary oscillations develop along the linearly increasing segment of the periodic as the locus \textbf{M} is approached. On the other hand, we find numerically that the parameter region for which we observe eventually two-dimensional dynamics extends from near \textbf{DZ} far into the $(\alpha,\beta)$-plane, as is suggested by the Poincar{\'e}-Bendixson-type result of Kennedy \cite{kennedy2019poincare}.
More specifically, the Floquet spectrum of the continued periodic orbits consists of two real Floquet multipliers with all others virtually indistinguishable from zero.

Overall, our findings may serve as a `health warning' to modelers who are faced with systems where state-dependency arises; concrete examples are DDE models for human balancing \cite{Insperger2011,ims2012}, for machining and milling \cite{IST07}, for laser systems \cite{MartinezLlinas2015}, as well as conceptual climate models of delayed action oscillator type \cite{Keane2019}. At the same time, our work demonstrates that state-dependent DDEs, including those with delay inside the state dependence, can be investigated effectively with advanced continuation tools, as implemented in the package DDE-BIFTOOL, in combination with analytical approaches including normal form analysis. In other words, researchers need not shy away from DDE models of this form when they arise in an application context, thus, allowing for the full investigation of the effects of different types of delayed feedback loops on observable behavior. Indeed, there will be quite a number of DDEs that can be studied in this spirit, and we hope that the work presented here will serve as a motivation to other researchers. This includes further work on saw-tooth shaped or other types of limits of periodic orbits, where numerical investigations may well stimulate further analysis, including into the effects of small additional delays in the state-dependent term.

\section*{Acknowledgments}

We thank Andrew Keane for helpful discussions and for performing some initial continuations. A.R.H.~is funded by the Natural Sciences and Research Council of Canada Discovery Grant RGPIN-2018-05062.
The research of B.K.~was supported by the Royal Society Te Ap\={a}rangi Marsden Fund grant \#19-UOA-223.
J.S.’s research is supported by the UK Engineering and Physical Sciences Research Council (EPSRC) grants EP/N023544/1 and EP/V04687X/1.



\providecommand{\href}[2]{#2}
\providecommand{\arxiv}[1]{\href{http://arxiv.org/abs/#1}{arXiv:#1}}
\providecommand{\url}[1]{\texttt{#1}}
\providecommand{\urlprefix}{URL }


\appendix
\section{Linear systems for the expansion coefficients of the center manifold}
\label{app:ex}

We provide here technical details of how to compute expansions of the ODE on the center manifold near the point \textbf{DZ}, given in Lemma~\ref{thm:cmf2}, up to any order. The basis $B$ of the subspace corresponding to the double eigenvalue $0$ of the linear operator $A$, given in \eqref{eq:derA}, is a row vector of length $2$ of the form
\begin{align*}
B:[-\tmax,0] &\to R^{1\times2} \\ 
\theta &\mapsto B(\theta)=[1,\theta]\mbox{.}
\end{align*}
Thus, $B$ maps from $\R^2$ into $C = C^0$. The spectral projection of an arbitrary function $u\in C$ onto $\lspan B$ is given via the functional
\begin{align*}
        \badj:C &\to \R^2\\
        u & \mapsto
	\begin{bmatrix}
		2/3\\ 2
	\end{bmatrix}u(0)+\int_{-1}^0
	\begin{bmatrix}
		2s+4/3\\ -2
	\end{bmatrix}u(s)\d s\mbox{.}
\end{align*}
Thus, $B\badj:C\to\lspan B\subset C$ is the spectral projection in $C$. The spectral projection of the operator $A$ corresponding to the right-hand side of the linearized DDE \eqref{lindde} with $\alpha=1$ (see \eqref{eq:derA}) onto the basis $B$ is
\begin{displaymath}
J_0=\badj A B=
\begin{bmatrix}
0&1\\0&0
\end{bmatrix}\in\R^{2\times 2}\mbox{.}
\end{displaymath}
As \cite{Hum-Dem-Mag-Uph-12,Sieber2017} laid out, all formal expansion coefficients of the center manifold $h:\R^2\times\R^2\to C$ are well defined to arbitrary order. Let us express the graph $h$ in its expansion coefficients (where we use $;$ for separating entries of column vectors)
\begin{displaymath}
h(x_c;p;q)(\theta)=Bx_c+\sum_{j=2}^m h^j(\theta)(x_c;p;q)^j+o(|(x_c;p;q)|^{m+1})
\end{displaymath}
with the unknown symmetric multilinear coefficients $h^j:[-\tmax,0]\mapsto \R^{1\times 4^j}$. We take into account the dependence of $h$ on the parameters $p$ and $q$, expanding $h$ in these parameters, too. The resulting ODE on the center manifold has the form
\begin{equation}
\label{expxdot}
\dot x_c
= J_0x_c+
\sum_{j=2}^mA^j_c
(x_c;p;q)^j+o(|(x_c;p;q)|^m)\mbox{,}
\end{equation}
where the $A^j_c\in\R^{2\times 4^j}$ are also unknown symmetric multilinear coefficients. The unknown coefficients $A^j_c$ and $h^j(\theta)$ are determined by comparing coefficients at order $j\geq2$ for the expansion of the invariance equation for $h$, \eqref{inv:de} and \eqref{inv:bc} below, and an orthogonality condition (\eqref{b:orth} below), given by
\begin{align}
	\label{inv:de}
	\partial_{x_c}h((x_c;p;q),\theta)\dot x_c&=\partial_\theta h((x_c;p;q),\theta)\mbox{,}\\
	\label{inv:bc}
	\partial_\theta h((x_c;p;q),0)&=f_{p,q}(h((x_c;p;q),\cdot))\mbox{,}\\
	\badj h^j(\theta)(x_c;p;q)^j&=0\mbox{.}
\label{b:orth}
\end{align}
The center manifold, extended by the parameters $p$ and $q$, is four-dimensional, with two trivial equations $\dot p=0$ and $\dot q=0$, such that we may extend $J_0$ as
\begin{displaymath}
J_{0,e}=\operatorname{diag}(J_0,0,0)\in\R^{4\times 4}\mbox{.}
\end{displaymath}
Thus, the unknowns at order $j$ are $h^j(\theta)e_{k_1}\cdots e_{k_j}\in\R$ and $A_c^je_{k_1}\cdots e_{k_j}\in \R^2$, where $e_\ell$ ($\ell=1,\ldots,4$) are the unit vectors in $\R^4$ and $(k_1,\ldots,k_j)$ is a non-decreasing element of 
${\mathcal J}=\{1,\ldots,4\}^j$ ($k_1\leq k_2\leq k_3\leq k_4$).

The matrix $J_0$ is upper triangular. Thus, putting the non-decreasing elements of $\{1,\ldots,4\}^j$ into lexicographic order, system \eqref{inv:de}--\eqref{b:orth} generates a sequence of $3$-dimensional linear systems of equations for each element $\kappa=(k_1,\ldots,k_j)$ of ${\mathcal J}$ and the corresponding pair of unknowns
$h^j_\kappa(\theta)=h^j(\theta)e_{k_1}\cdots e_{k_j}\in\R$ and
$A_{c,\kappa}^j=A_c^je_{k_1}\cdots e_{k_j}\in \R^2$.

Using the expansion of $h$ and $f$, equation~\eqref{inv:de} results in a linear ODE of the form for each order $j$ and coefficient $\kappa$
\begin{equation}\label{hj:type}
h^j_\kappa(\theta)'=\Lambda_\kappa\cdot h^j_\kappa(\theta)+B(\theta)A_{c,\kappa}+r_{h,\kappa}(\theta)\mbox{,}
\end{equation}
where $\Lambda_\kappa=\sum_{\ell\in\kappa}\lambda_\ell$ is the sum of eigenvalues of $J_{0,e}$ corresponding to index set $\kappa$ (thus, $\Lambda_\kappa=0$ for all $\kappa$ in our case). The inhomogeneity $r_{h,\kappa}(\theta)$ comes from previously determined lower-order terms $h^\ell$ with $\ell<j$ and elements $h^j_{\kappa'}$, where $\kappa'<\kappa$ in the lexicographic order.

Equation~\eqref{inv:bc} and its expansion provides a scalar equation for each order $j$ and coefficient $\kappa$, involving the unknowns $h^j_\kappa(0)$ and $A^j_{c,\kappa}$. Its coefficients in our concrete example are $-L[\theta\mapsto1]=0$ for $h^j_\kappa(0)$ and $B'(0)-L[\theta\mapsto \int_0^\theta B(s)\d s=(-1,3/2)$ for $A^j_{c,\kappa}$. The terms involving higher-order derivatives of $f_{(p,q)}$ all depend on previously computed lower order terms $h^\ell(\theta)$ ($\ell<j$) only. Since $f_{(p,q)}$ is only mildly differentiable, their evaluation requires differentiation of $h^\ell(\theta)$ with respect to $\theta$. However, since the $h^\ell$ are all solutions of linear ODEs of the form \eqref{hj:type}, the arguments of $f_{(p,q)}$ are in $C^\infty$, making the expansion of $f_{(p,q)}(h((x_c;p;q),\cdot))$ possible to arbitrary order.

The remaining equations are enforced by the orthogonality condition \eqref{b:orth} on $h^j_\kappa(\theta)$. This is in contrast to the approach taken in \cite{Kuz-04-Book} and \cite{BJK20}, where the coefficients $A^j_c$ are forced to be in an a-priori known normal form for a range of standard bifurcations. Since there is no well-known normal form established for the degenerate point \textbf{DZGH}, we keep the terms $h^j(\theta)(x_c;p;q)^j$ orthogonal to the linear spectral projection $\badj$. In our concrete example, the coefficients for $h^j_\kappa(0)$ are $\badj[\theta\mapsto1]=(1;0)$, and, for $A^j_{c,\kappa}$, $\badj[\theta\mapsto\int_0^\theta B(s)\d s]=(0,-1/36;1,-1/3)$, such that the overall coefficient matrix for the unknowns $(h^j_\kappa(0);A^j_{c,\kappa})$ at all orders $j>1$ and for all index sets $\kappa$ is the invertible matrix
\begin{displaymath}
M=
\begin{bmatrix}
0&-1&\phantom{-}3/2\\ 1&\phantom{-}0&-1/36\\0&\phantom{-}1&-1/3
\end{bmatrix}\mbox{.}
\end{displaymath}
The complete manifold coefficient $\theta\mapsto h^j_\kappa(\theta)$ can then be obtained by solving \eqref{hj:type} from the initial condition $h^j_\kappa(0)$ with the known coefficients $A_{c,\kappa}$ an the inhomogeneity $r_{h,\kappa}(\theta)$ from lower-order terms. 

\section{Contents of supplementary material}
\label{app:supp}

The file available at \url{http://auckland.figshare.com/articles/online_resource/hkrs_sdDDE_cmf_supplement_zip/16735474} provides general Matlab and Octave (sympy) compatible symbolic algebra drivers for the expansion of center manifolds in DDEs with state-dependent discrete delays. It implements the expansion described in Appendix~\ref{app:ex} and demonstrates it up to truncation order $5$ for the specific example of the point \textbf{DZGH} of \eqref{eq:sdDDE}. The symbolic expressions that are generated as output can be converted to right-hand sides for use in numerical bifurcation analysis packages.


\end{document}